\documentclass[11pt]{article}
\usepackage{amsmath,amssymb,amsthm,graphicx,pstricks,xypic,comment,fancyhdr}

\title{Bending Fuchsian representations of fundamental groups of
  cusped surfaces in PU(2,1).}  \author{Pierre WILL\\Institut
  Fourier\\100 rue des Maths\\38402 St Martin d'H\`eres\\
  France\\\texttt{pierre.will@ujf-grenoble.fr}}

\begin{document}
\maketitle
\def\og{\leavevmode\raise.3ex\hbox{$\scriptscriptstyle\langle\!\langle$~}}
\def\fg{\leavevmode\raise.3ex\hbox{~$\!\scriptscriptstyle\,\rangle\!\rangle$}}
\newcommand{\la}{\langle}
\newcommand{\ra}{\rangle}
\newcommand{\HdC}{\bold H^{2}_{\mathbb C}}
\newcommand{\HdR}{\bold H^{2}_{\mathbb R}}
\newcommand{\HtR}{\bold H^{3}_{\mathbb R}}
\newcommand{\HuC}{\bold H^{1}_{\mathbb C}}
\newcommand{\HnC}{\bold H^{n}_{\mathbb C}}
\newcommand{\Cdu}{\mathbb C^{2,1}}
\newcommand{\Ct}{\mathbb C^{3}}
\newcommand{\C}{\mathbb C}
\newcommand{\Mn}{\mbox{M}_n\left(\mathbb C\right)}
\newcommand{\R}{\mathbb R}
\newcommand{\A}{\mathbb A}
\newcommand{\cC}{\mathfrak C}
\newcommand{\cCp}{\mathfrak C_{\mbox{p}}}
\newcommand{\cT}{\mathfrak T}
\newcommand{\rR}{\mathfrak R}
\newcommand{\hH}{\mathfrak H}
\newcommand{\n}{\noindent}
\newcommand{\p}{{\bf p}}
\newcommand{\Fa}{\mathcal{F}_\infty}
\newcommand{\tZ}{\texttt{Z}}
\newcommand{\lox}{\mbox{\scriptsize{lox}}}
\newcommand{\Td}{T^{2}_{\left(1,1\right)}}
\newcommand{\T}{T_{\left(1,1\right)}}
\newcommand{\PSL}{\mbox{\rm PSL(2,$\R$)}}
\newcommand{\Pu}{\mbox{\rm PU(2,1)}}
\newcommand{\X}{\mbox{\textbf{X}}}
\newcommand{\tr}{\mbox{\rm tr}}
\newcommand{\B}{\mathcal{B}}
\newcommand{\bp}{{\bf p}}
\newcommand{\bm}{{\bf m}}
\newcommand{\bc}{{\bf c}}
\newcommand{\bn}{{\bf n}}
\newcommand{\bq}{{\bf q}}
\newcommand{\td}{{\tt d}}
\newcommand{\tD}{{\tt D}}
\newcommand{\reppudu}{Rep$_{\pi_1(\Sigma),{\rm PU(2,1)}}$}
\newcommand{\repg}{Rep$_{\pi_1(\Sigma),G}$ }
\newcommand{\reppsldr}{Rep$_{\pi_1(\Sigma),{\rm PSL}(2,\R)}$ }
\renewcommand{\arg}{\mbox{arg}}
\renewcommand{\Re}{\mbox{\rm Re}\,}
\renewcommand{\Im}{\mbox{\rm Im}\,}
\renewcommand{\leq}{\leqslant}
\renewcommand{\geq}{\geqslant}  
\newcommand{\tol}{\textbf{tol}}
\newtheorem{theo}{Theorem}
\newtheorem*{theo*}{Theorem}
\newtheorem{lem}{Lemma}
\newtheorem{coro}{Corollary}
\newtheorem{prop}{Proposition}
\theoremstyle{remark}
\newtheorem{rem}{Remark}
\newtheorem{ex}{Example}
\theoremstyle{definition}
\newtheorem{defi}{Definition}

\begin{abstract}
  We describe a new family of representations of $\pi_1(\Sigma)$ in
  PU(2,1), where $\Sigma$ is a hyperbolic Riemann surface with at
  least one deleted point. This family is obtained by a bending
  process associated to an ideal triangulation of $\Sigma$. We give an
  explicit description of this family by describing a coordinates
  system in the spirit of shear coordinates
  on the Teichm\"uller space. We identify within this family new
  examples of discrete, faithful and type-preserving representations
  of $\pi_1(\Sigma)$. In turn, we obtain a 1-parameter family of
  embeddings of the Teichm\"uller space of $\Sigma$ in the
  PU(2,1)-representation variety of $\pi_1(\Sigma)$.  These results
  generalise to arbitrary $\Sigma$ the results obtained in
  \cite{Wi2} for the 1-punctured torus.
\end{abstract}

{\small Key words: Complex hyperbolic geometry, representations of
  surface groups, Teichm\"uller space, deformation.}

{\small AMS classification 51M10, 32M15, 22E40}

\tableofcontents
\section{Introduction}
Let $\Sigma$ be an oriented surface with negative Euler
characteristic. Describing the representation variety of the fundamental group
$\pi_1(\Sigma)$ in a given Lie group $G$ has been a major problem
during the last two decades. The central object in  this field is the 
character variety

\begin{equation}\mbox{\repg}={\rm Hom}(\pi_1(\Sigma),G)//G\label{rep}.
\end{equation}

This problem finds its source in the study of the Teichm\"uller space
of $\Sigma$, which classifies hyperbolic metrics or complex structures
on $\Sigma$. Riemann's uniformization theorem implies that the
Teichm\"uller space of $\Sigma$ may be seen as the subset of \reppsldr
consisting of conjugacy classes of discrete, faithful and
type-preserving representations.  In the case where $\Sigma$ is closed
without boundary, Goldman classified in \cite{Go5} the connected
components of \reppsldr using the Euler number of a representation. It
turns out that the extremal values of the Euler number correspond to
two connected components of \reppsldr which are copies of the
Teichm\"uller space of $\Sigma$. Since then, the question of
classifying the connected components of \repg, and understanding the
situation of dicrete and faithful representations has been addressed
for many Lie groups (see for instance
\cite{BuIoWiLab,BuIoWi,FockGon2,Hit,Lab}). This led to what is
sometimes called \textit{higher Teichm\"uller theory}.  The goal of
this work is to present an explicit family of geometrically
well-understood representations of the fundamental group of a cusped
surface in PU(2,1), the group of holomorphic isometries of the complex
hyperbolic plane $\HdC$.

In the specific context of PU($n$,1) and complex hyperbolic geometry,
the study of representations of surface groups has been initiated by
Goldman and Toledo among others (\cite{GM,Go7,Tol}).  Since then, the
question of the classification of connected components of \reppudu has
been answered thanks to Toledo's invariant (see \cite{Tol,Xia}). This
invariant is defined as the integral over $\Sigma$ of the pull-back of
the K\"ahler form on $\HdC$ by a $\rho$-equivariant embedding of
$\tilde \Sigma$ into $\HdC$. We will denote it by $\tol(\rho)$. If
$\Sigma$ is compact without boundary, the Toledo invariant enjoys the
following properties.

\begin{itemize}
\item $\tol$ is a continuous function of $\rho$.
\item $\forall\rho\in\,$\reppudu, $\tol(\rho)\in 2/3\mathbb{Z}$.
\item $\tol$ satisfies a Milnor-Wood inequality: $|\tol(\rho)|\leq
  -4\pi\chi(\Sigma)$.
\end{itemize}

Xia proved in \cite{Xia} that two representations having the same
Toledo invariant lie in the same connected component of \reppudu, and
thus the Toledo invariant classifies connected components of \reppudu.
 Extremal values of the Toledo invariant characterizes
those representations preserving a complex line on which they
act properly discontinuously (see \cite{GM,Tol}). These results gave rise to
considerable generalisations, from the context of complex hyperbolic
space to the wider frame of Hermitian symmetric spaces (see for
instance \cite{BuIoWiLab,BuIoWi,KoMa,KoMa2}).

The question of discreteness of representations of surface groups in
PU($n$,1) is still far from being given a complete answer. It is known
for instance that contrary to the case of PSL(2,$\R$), discrete and
faithful representations are not contained in specific components of
Rep$_{\Sigma,{\rm PU(2,1)}}$, as shown for instance in \cite{GKL} or
\cite{AGG1,AGG2,Gaye}.

This work is concerned with representations of fundamental groups of
cusped surfaces. In this case, the Toledo invariant is defined for
type-preserving representation (that is, representations mapping
classes of loops around punctures to parabolics), but it does not
classify the topological components of \reppudu, as shown in
\cite{GuP2}.  The question of the discreteness for such groups has
been addressed in several works, see for example
\cite{FK1,FK,GP,GuP,GuP2,S1,S2,S3,Wi2}.

A classical way to produce non-trivial examples of representations of
surface groups in PU(2,1) is to start with a representation $\rho_0$
preserving a totally geodesic subspace $V$ of $\HdC$ and to deform it.
There are a priori two ways to do so in our case since there are two
kinds of maximal totally geodesic subspaces in $\HdC$, namely
\textit{complex lines} and totally geodesic Lagrangian planes, or
\textit{real planes}. Complex lines are embeddings of $\HuC$ with
sectional curvature $-1$, and real planes are embeddings of $\HdR$
with sectional curvature $-1/4$. Their respective stabilisers are
P(U(1,1)$\times$U(1)) and PO(2,1) (see for instance \cite{Go}). A
discrete and faithful representation preserving a complex line (resp.
a real plane) is called $\C$-\textit{Fuchsian}
(resp. $\R$\textit{-Fuchsian}). The rigidity results for $\C$-Fuchsian
representations in the case of compact surfaces are no longer true for
$\R$-Fuchsian representation see \cite{Gui,PP1}.\\ 

\n The purpose of this work is twofold.

\begin{enumerate}
\item We first describe a family of representations obtained by a
  bending process. They arise as holonomies of equivariant mappings
  from the Farey set of a surface $\Sigma$ to the boundary of
  $\HdC$. Roughly speaking, we are bending along the edges of an ideal
  triangulation of $\Sigma$, and the case were the representation is
  $\R$-Fuchsian corresponds to vanishing bending angles (Theorem
  \ref{realbend}).
 
\item We identify within this family of examples a subfamily of
  discrete, faithful, and type-preserving representations (Theorem
  \ref{theodiscret}). The proof of discreteness is done by showing that
  the action of $\rho(\pi_1)$ on $\HdC$ is discontinuous. We obtain in
  turn a 1-parameter family of embeddings of the Teichm\"uller space
  of $\Sigma$ into \reppudu containing only classes of discrete,
  faithful representations with unipotent boundary holonomy (Theorem
  \ref{theoembed}).
\end{enumerate}

\n In the context of $\HdC$, deformations by bending were first
described by Apanasov in \cite{Apa}. More recently, in \cite{Plat},
Platis has described a complex hyperbolic version of Thurston's
quakebending deformations for deformations of $\R$-Fuchsian
representations of groups in the case of closed surfaces without
boundary. If $\rho_0$ is an $\R$-Fuchsian representation of
$\pi_1(\Sigma)$ and $\Lambda$ is a finite geodesic lamination with a
complex transverse measure $\mu$, Platis shows that there exists
$\epsilon>0$ such that any quakebend deformation $\rho_{t\mu}$ of
$\rho_0$ is complex hyperbolic quasi-Fuchsian for all $t<\epsilon$.
The proof of discreteness in \cite{Plat} rests on the main result in
\cite{PP1} where the proof of discreteness is done by building a fundamental 
domain. \\

\n In order to sum up our work, let us introduce a little notation.\\

Throughout this work, we denote by $\Sigma$ a oriented surface of
genus $g$ with $n>0$ deleted points, which we denote by $x_1,\cdots
x_n$. We assume that $\Sigma$ has negative Euler characteristic,
that is, $2-2g-n<0$. We denote by $\pi_1(\Sigma)$ the
fundamental group of $\Sigma$. It admits the following presentation

$$\pi_1(\Sigma)\sim\la a_1,b_2,\cdots a_g,b_g,c_1\cdots c_n\lvert \prod [a_i,b_i]\prod c_j=1\ra,$$

\n where the $c_j$'s are the homotopy classes of simple loops
enclosing the punctures $x_j$ of $\Sigma$. The universal cover of
$\Sigma$ is an open disc $\tilde\Sigma$, with a $\pi_1$-invariant
family of points on its boundary corresponding to the deleted points
of $\Sigma$. This set of boundary points is called the \textit{Farey
  set} of $\Sigma$, and denoted by $\mathcal{F}_\infty$. If one fixes
a finite area hyperbolic structure on $\Sigma$, the Farey set consists
of the fixed points of parabolic elements representing homotopy
classes of loops around punctures. In particular, if the holonomy of
this hyperbolic structure has its image in PSL(2,$\mathbb{Z}$), one
recovers this way the usual Farey set $\mathbb{Q}\cup\infty$.

The idea of shear coordinates on the Teichm\"uller space of a cusped
surface goes back to the eighties with Thurston \cite{Thu2}, Bonahon
\cite{Bonah} and Penner \cite{PenR}. The principle is the
following. In order to build a hyperbolic structure on a cusped
surface $\Sigma$, it suffices to glue together ideal triangles in the
upper half plane. Since there is a unique ideal triangle up to the
action of \PSL, the only gluing invariant is the cross ratio of the
four boundary points associated to a pair of adjacent triangles, which
measures the shearing of the two triangles. It is therefore very
natural to parametrize structures by decorating triangulations, using
(positive) real numbers that should be interpreted as cross ratios.
In particular, Penner defined the decorated Teichm\"uller space of a
punctured surface in \cite{PenR}, and gave a description of it in this
way. Later, these ideas were successfully exploited and generalised by
Fock and Goncharov to study representations of cusped surfaces in real
split Lie groups (see \cite{FockGon2,FockGonA,FockGon3}). It is
possible to give an efficient combinatorial description of a
representant of the class of representation of $\pi_1(\Sigma)$
associated to a given decoration (see for instance \cite{FockGonA}),
which we are going to adapt to our context (section \ref{explicit}).

As shown in \cite{MW}, it is possible to describe a similar system of
explicit coordinates on an open subset of \reppudu which contains all
the classes of discrete and faithful representations. However,
identifying those classes of representations that are indeed discrete
remains out of reach. Therefore, we restrict ourselves to a
family of representations obtained by making an additional geometric
assumption.  More precisely, our first goal is to classify what we
call \textit{$T$-bent realizations} of $\mathcal{F}_\infty$ in $\HdC$,
where $T$ is an ideal triangulation of $\Sigma$, that is pairs
$(\phi,\rho)$, where

\begin{itemize}
\item $\rho$ is a representation $\pi_1(\Sigma)\longrightarrow$ Isom($\HdC$),
\item $\phi:\mathcal{F}_\infty\longrightarrow \partial\HdC$ is a
  $(\pi_1,\rho)$-equivariant mapping,
\item for any face $\Delta$ of $\hat T$ with vertices $a,b,c\in
  \mathcal{F}_\infty$, the three points $\phi(a)$, $\phi(b)$ and
  $\phi(c)$ form a real ideal triangle of $\HdC$, that is, they
  belong to the boundary of a real plane of $\HdC$.
\end{itemize}

It is natural to use this formalism because of the following remark.
Consider a hyperbolic structure on $\Sigma$ with holonomy $\Gamma$ a
subgroup of PSL(2,$\mathbb{Z}$). Each point of the Farey set
$\mathbb{Q}\cup\infty$ is the fixed point of a unique primitive
parabolic element of $\Gamma$, corresponding to a peripheral loop. Now
taking $\Gamma$ as a Fuchsian model for $\pi_1(\Sigma)$, any other
finite area hyperbolic structure on $\Sigma$ gives rise to a discrete,
faithful and type preserving representation
$\rho:\Gamma\longrightarrow\PSL$. One associates to $\rho$ an
equivariant mapping $\phi:\Fa\longrightarrow \partial\HuC$ by sending
any point $m$ in the Farey set corresponding to the fixed point of the
primitive parabolic $p$ to the fixed point of $\rho(p)$ (see also
section \ref{embed}).

In order to parametrize these bent realizations, we need a gluing
invariant for pairs of real ideal triangles in $\HdC$. In the context
of PSL(2,$\R$), this invariant is the cross-ratio: the unique
invariant of a pair of ideal triangles of $\HuC$ sharing an edge. In
the complex hyperbolic context, a pair of real ideal triangles sharing
an edge in $\HdC$ is determined up to holomorphic isometry by a single
complex number $\tZ(\tau_1,\tau_2)\in\C\setminus \lbrace
-1,0\rbrace$. This $\tZ$-invariant is similar to the Kor\'anyi-Reimann
cross-ratio on the Heisenberg group (see \cite{Go,KR,Wi4} and remark
\ref{simKR} in section \ref{invZ}). Note that $\tZ$ is an invariant of
ordered pairs of real ideal triangles, in the sense that

\begin{equation}\label{z-orient}
\tZ(\tau_2,\tau_1)=\overline{\tZ(\tau_1,\tau_2)}.
\end{equation}

This invariant was first used by Falbel in \cite{F} to glue ideal
tetrahedra in $\HdC$. Falbel needs two such parameters to describe the
isometry class of an ideal tetrahedron.  We only need one since we only
consider pairs of real ideal triangles, which correspond in his
terminology to the particular case of symmetric tetrahedra (see
section 4.3 of \cite{F}).

The modulus of $\tZ$ is similar to the cross-ratio in $\HuC$ : it
measures the shearing between two real ideal triangles. Its argument
is the \textit{bending parameter}, which can be seen as the measure of
an angle between real planes. In particular, if $\tZ$ is real, the two
adjacent real ideal triangles are contained in a common real plane.
We will therefore call a \textit{bending decoration} of $T$ any
application $\tD: e(T)\longrightarrow \C\setminus\lbrace -1,0\rbrace$
(the two cases where $z=0$ or $z=-1$ correspond to degenerate
triangles).

As in the case of PSL(2,$\R$), it is possible to associate $T$-bent
realizations to bending decorations and we obtain an explicit
expression for the images of classes of loops by $\rho$. We do this in
the same spirit as in \cite{FockGonA}. Once an ideal triangulation of
$\Sigma$ is chosen, any element $\gamma$ of $\pi_1(\Sigma)$ can be
represented by a sequence of edges of the modified dual graph (see
Definition \ref{dual}). Using the decoration, we associate to each
edge of this graph an elementary isometry, and the image of
$\gamma$ by $\rho$ is the product of the corresponding elementary
isometry. These elementary isometries are either elliptic elements of
order 3, or antiholomorphic involutions. The latter appear as the
unique isometries exchanging two real ideal triangles with a common
edge. They have to be antiholomorphic because of relation
(\ref{z-orient}). As a consequence the image of a homotopy class by
$\rho$ is not always holomorphic. This is why the representation
$\rho$ is taken in Isom($\HdC$) rather than in PU(2,1), which is the
index two subgroup of Isom($\HdC$) containing holomorphic
isometries. However, for some special triangulations, the image of the
representation is in fact contained in PU(2,1). Namely, we show in
section \ref{when}, that the representation $\rho$ associated to a
$T$-bent realization of $\mathcal{F}_\infty$ is holomorphic if and
only if $T$ is \textit{bipartite}, that is if its dual graph is
bipartite . Now, any cusped surface $\Sigma$ admits a bipartite ideal
triangulation (Proposition \ref{existbipar}).  This bending
process produces thus representations of $\pi_1(\Sigma)$ in PU(2,1)
for any non compact $\Sigma$. Let us denote by $\mathcal{BD}_T$ the
set of bending decorations of an ideal triangulation $T$ and by
$\mathcal{BR}^*_T$ the quotient of $\mathcal{BD}_T$ by the action of
complex conjugation. The first result of our work is the following.\\

\begin{theo}[\textbf{Bending Theorem}]\label{realbend}
  There is a bijection between $\mathcal{BD}^*_T$ and
  $\mathcal{BR}_T$.
\end{theo}

\n More precisely, we associate to any bending decoration a 
pair $(r_1,r_2)$ of PU(2,1)-classes of bent realizations of
$\mathcal{F}_\infty$, which represent the same Isom($\HdC$)-class of
realization. The complex conjugation of bending representations
corresponds to the permutation $(r_1,r_2)\longrightarrow (r_2,r_1)$.
This result is proved in section \ref{proof2}. 

After having classified $T$-bent realizations, we focus on a special
kind: those $T$-bent realizations corresponding to \textit{regular}
bending decorations, that is decorations of the form $\tD=\td
e^{i\theta}$ , where $\td$ is a positive decoration of $T$, and
$\theta\in]-\pi,\pi]$ is a fixed real number. When $\theta=0$, we
obtain $T$-bent realizations where all the images of the points of
$\mathcal{F}_\infty$ are contained in a real plane. The corresponding
representations are $\R$-Fuchsian.\\

\n If $c$ is a class of peripheral loop on $\Sigma$ the parabolicity
of $\rho(c)$ is simply expressed in terms of the decoration. If $x$ is
the puncture of $\Sigma$ surrounded by $c$, and $e_1,\cdots,e_n$ are
the edges of $T$ adjacent to $x$, the isometry $\rho(c)$ is parabolic
if and only if the product $\prod \tD(e_i)$ has modulus 1. Such a
decoration $\tD$ is said to be \textit{balanced at $x$} (see section
\ref{holes}).

\n We now state the main result of our work.\\

\begin{theo}[\textbf{Discreteness Theorem}]\label{theodiscret}
  Let $T$ be a bipartite ideal triangulation of $\Sigma$,
  $\theta\in]-\pi,\pi[$ be a real number and $\tD$ be a regular
  bending decoration of $T$ with angular part equal to $\theta$. Let
  $\rho$ be a representative of the (unique) Isom($\HdC$)-class of
  representations associated to $\tD$. Then
\begin{enumerate}
\item For any index $i$, $\rho(c_i)$ is parabolic if and only if $\tD$
  is balanced at $x_i$.
\item The representation $\rho$ does not preserve any totally geodesic
  subspace of $\HdC$, unless $\theta=0$, in which case it is
  $\R$-Fuchsian.
\item As long as $\theta\in[-\pi/2,\pi/2]$, the representation $\rho$
  is discrete and faithful.
\end{enumerate}
\end{theo}

 \n The bending and discreteness theorems are generalisations to the
case of any punctured surface of results obtained in \cite{Wi2}
in the case of the punctured torus.

\n If we restrict this result to those decorations $\tD$ which are
balanced at every puncture, we obtain a 1-parameter family of
embeddings of the Teichm\"uller space of $\Sigma$ in \reppudu
$\pi_1(\Sigma)$. The images of these embedding contain
only classes of discrete, faithful and type-preserving
representations, parametrised by $\theta\in[-\pi/2,\pi/2]$. Note that
the family of embeddings obtained depends on the initial choice
of the triangulation.

\n The proof of Theorem \ref{theodiscret} goes as follows. From a
bending decoration, we construct a family of real ideal triangles on
which $\rho(\pi_1)$ acts because of the equivariance condition.  Under
the hypotheses of Theorem \ref{theodiscret},we are able to define for
each pair $(\tau,\tau')$ of real ideal triangles a canonical real
hypersurface of $\HdC$ called the \textit{splitting surface} of $\tau$
and $\tau'$ and denoted by Spl$(\tau,\tau')$. The main steps of the
proof consist in proving the following properties for these
hypersurfaces.

\begin{enumerate}
\item The splitting surface Spl$(\tau,\tau')$ separates $\HdC$ in two
  connected components, each of which contains one of $\tau$ and
  $\tau'$ (Proposition \ref{splitplits}).
\item If $\tau$ is surrounded by three real ideal triangles
  $\tau_1$, $\tau_2$ and $\tau_3$, the three splitting surfaces
  Spl$(\tau,\tau_i)$ are mutually disjoint (Theorem
  \ref{crucial}).
\end{enumerate}

It is a direct consequence of these facts that all the constructed
triangles are disjoint, and $\rho(\pi_1(\Sigma))$ acts discontinuously
on their union, and is therefore discrete. The union of all the
triangles constructed is a piecewise totally geodesic disc.

Splitting surfaces are examples of spinal $\R$-surfaces (see section
\ref{defispinal}), which are the inverse images of geodesics by the
orthogonal projection on real planes. This terminology refers to
Mostow's spinal surfaces, defined in \cite{Most}, which are the
inverse images of geodesic by the orthogonal projection on a complex
line (spinal surfaces are often called \textit{bisectors}, see
\cite{Go}). Spinal $\R$-surfaces appeared first in \cite{Wi2} under
the name of $\R$-balls. They were generalised and used by Parker and
Platis in \cite{PP1} (see also \cite{PP2}) under the name of
\textit{packs}. In their terminology, spinal $\R$-surfaces are
\textit{flat packs}. In particular, the characterisation of spinal
$\R$-surfaces given in the Lemma \ref{orbit} is similar to their
definition of packs.

If $\Sigma$ has genus $g$ and $n$ punctures, any ideal triangulation
of $\Sigma$ has $6g-6+3n$ edges. The set of conjugacy classes discrete
and faithul representations of $\pi_1(\Sigma)$ (resp. the Teichmüller
space of $\Sigma$) has (real) dimension $6g-6+3n$ (resp. $6g-6+2n$).
The representation variety \reppudu has real dimension $16g-16+8n$,
and its subset containing the classes of type preserving
representations has real dimension $16g-16+7n$. If $T$ is an ideal
triangulation of $\Sigma$, $\mathcal{BR}_T$ and $\mathcal{BR}^*_T$
have real dimension $12g-12+6n$, and may be seen respectively as
$\left(\C\setminus\lbrace -1,0\rbrace\right)^{6g-6+3n}$ and its
quotient by the action of the complex conjugation.  The real
dimension of the family of the classes of discrete and faithful
representations obtained by examining regular bending decorations of
$T$ is $6g-6+3n+1$ and falls to $6g-6+2n+1$ if we add the condition of
type-preservation.  The Toledo invariant ${\bf tol}$ (see
\cite{GuP2,Tol}) of a representation is defined for representations of
fundamental groups of compact surfaces, and for type preserving
representations of surfaces with deleted points. All type-preserving
representations we obtain here have vanishing Toledo invariant. This
follows from the fact that the representations are constructed from
families of real ideal triangles (see section \ref{rem}).

Our work is organised as follows. We provide in section \ref{HdCisom}
the necessary background about the complex hyperbolic plane and its
isometries. The invariant of a pair of real ideal triangles is
described in section \ref{RIT}. Section \ref{sectionbending} is
devoted to the proof of Theorem \ref{realbend}.We provide in
\ref{explicit} the explicit form of the corresponding representations.
The characterisation of bent realization giving representations in
PU(2,1) in terms of bipartite triangulations is given in \ref{when},
and we study the holonomy of loops around deleted points in
\ref{holes}. We turn then to the proof of Theorem
\ref{sectiondiscret}. Spinal $\R$-surfaces and splitting surfaces
surfaces are defined and studied in \ref{defispinal}, and we prove the
discreteness part of Theorem \ref{theodiscret} in \ref{fin}. Section
\ref{rem} is devoted to some remarks and comments. In particular, give
a quick presentation of the similar classical construction for \PSL,
and prove Theorem \ref{theoembed}.We draw the connection between our work
and the previously known families of examples studied in
\cite{FK1,GuP2,Wi2}.

\paragraph{Acknowledgements} I would like to thank Nicolas Bergeron,
Julien March\'e and Anne Parreau for fruitful discussions, and John
Parker for a useful hint about Proposition \ref{existbipar}. Part of
this work was carried out during a stay at the Max Planck Institut
f\"ur Mathematik in Bonn, and I would like to thank the institution
for the wonderful working conditions I had there. I thank the referees
for their help to improve this work. Last but not least,
I thank Elisha Falbel for his constant interest and support.

\section{The complex hyperbolic 2-space\label{HdCisom}}
We refer the reader to \cite{CG,Go} for more precise information and
more references about the material exposed in this section.

\subsection{$\HdC$ and its isometries}

Let $\Cdu$ denote the vector space $\C^3$ equipped with the Hermitian
form of signature (2,1) given by the matrix
\begin{equation}\label{J}
J=\begin{bmatrix}0 & 0 & 1\\0 & 1 & 0\\1 & 0 & 0\end{bmatrix}.
\end{equation}
The hermitian product of two vectors $X$ and $Y$ is given by $\la
X,Y\ra= X^TJ\bar Y$, where $X^T$ denotes the transposed of $X$. We
denote by $V^-$ (resp. $V^0$) the negative (resp. null) cone
associated to the hermitian form, and by $P$ the projectivisation $P :
\Cdu\longrightarrow \C P^2$.

In the rest of the paper, whenever $m$ is a point in $\C P^2$, we will
denote by $\textbf{m}$ a lift of it to $\Cdu$.
\begin{defi}
  The complex hyperbolic 2-space $\HdC$ is the projectivisation of
  $V^-$ equipped with the distance function $d$ given by
\begin{equation}\cosh^2{d\left(\dfrac{m,n}{2}\right)}=\dfrac{\la\bm,\bn\ra\la\bn,\bm\ra}
  {\la\bm,\bm\ra\la\bn,\bn\ra}\quad \forall (m,n)\in P(V^�)^2\label{dist}
\end{equation}
\end{defi}

\begin{prop}\label{isomHdC}
  The isometry group of $\HdC$ is generated by PU(2,1), the
  projective unitary group associated to $J$ and the complex
  conjugation.
\end{prop}

The group PU(2,1) is the group of holomorphic isometries of $\HdC$,
and is the identity component of Isom($\HdC$). The other component
contains the antiholomorphic isometries, all of which may be written
in the form $\phi\circ\sigma$, where $\phi$ is a holomorphic isometry
and $\sigma$ is the complex conjugation.

\paragraph{Horospherical coordinates}
The complex hyperbolic 2-space is biholomorphic to the unit ball of
$\C^2$, and its boundary is diffeomorphic to the 3-sphere $S^3$. The
projective model of $\HdC$ associated to the matrix $J$ given by
(\ref{J}) is often referred to as the \textit{Siegel model} of $\HdC$.
In this model, any point $m$ of $\HdC$ admits a unique lift to $\C^3$
given by

\begin{equation}\bm=\begin{bmatrix}-|z|^2-u+it\\z\sqrt{2}\\1\end{bmatrix}\mbox{, with } 
z\in\C,\,t\in\R\mbox{ and }u>0.\label{lifthoro}\end{equation}

The boundary of $\HdC$ corresponds to those vectors for which $u$
vanishes together with the point of $\C\bf P^2$ corresponding to the
vector $\begin{bmatrix} 1 & 0 & 0\end{bmatrix}^T$.  It may this be seen as the
the one point compactification of $\R^3$. The triple $(z,t,u)$ given
by (\ref{lifthoro}) is called the \textit{horospherical coordinates}
of $m$ (the hypersurfaces $\{u=u_0\}$ are the horospheres centred at
the point $\infty$ of $\partial\HdC$, which corresponds to the vector
$\begin{bmatrix} 1 & 0 & 0\end{bmatrix}^T$).

The boundary of $\HdC$ has naturally the structure of the
3-dimensional Heisenberg group, seen as the maximal unipotent subgroup
of PU(2,1) fixing $\infty$. We will call \textit{Heisenberg
  coordinates} of the boundary point with horospherical coordinates
$(z,t,0)$ the pair $[z,t]$. In these coordinates the group structure
is given by
 
 $$[z,t]\cdot [w,s] = [z+w,s+t+2\Im (z\bar w)]$$

 \paragraph{The ball model of $\HdC$} The same construction can be done
 using a different hermitian form of signature $(2,1)$ on $\C^3$.
 Using the special form associated to the matrix
 $J_0=\mbox{diag}(1,1,-1)$, we would obtain the so-called \textit{ball
   model} of the complex hyperbolic 2-space, which lead to a
 description of $\HdC$ as the unit ball $\C^2$.

\subsection{Totally geodesic subspaces\label{total}}
The maximal totally geodesic subspaces of $\HdC$ have real dimension
2. There are two types of such subspaces: the \textit{complex lines},
and the \textit{real planes}.
\paragraph{The complex lines.} These subspaces are the images under
projectivisation of those complex planes of $\C^3$ intersecting the
negative cone $V^-$. The standard example is the subset $C_0$ of
$\HdC$ containing points of horospherical coordinates $(0,t,u)$ with
$t\in\R$ and $u>0$.  This is an embedded copy of the usual Poincar\'e
upper half-plane. We will refer to this particular complex line as
$\HuC\subset\HdC$. All the other complex lines are the images of
$\HuC$ by an element of PU(2,1). Note that any complex line $C$ is
fixed pointwise by a unique holomorphic involutive isometry, called
the \textit{complex symmetry about $C$}.

\paragraph{The real planes.}  These subspaces are the images of the
Lagrangian vector subspaces of $\Cdu$ under projectivisation. The
standard example is the subset containing points of horospherical
coordinates $(x,0,u)$ with $x\in\R$ and $u>0$. The image of the mapping

\begin{equation} \label{embedhalfsup}
x+iu \longmapsto (x,0,u)
\end{equation}

\n is again an embedded copy of the usual Poincar\'e upper half-plane and
we will refer to this particular real plane as $\HdR\subset\HdC$. 
All other real planes are images of the standard one by an element of
PU(2,1).  There is also a unique involution fixing pointwise a real
plane $R$ which is called the \textit{real symmetry about $R$}. It is
antiholomorphic, and, in the case of $\HdR$, is complex
conjugation. If $\sigma$ is a real symmetry, we will call the real
plane which is its fixed point set its \textit{mirror}.

\begin{rem}
  In the ball model, the standard complex line is the first axis of
  coordinates $\lbrace (z,0),|z|<1\rbrace$. The standard real plane
  $\HdR$ is  the set of points with real coordinates $\lbrace
  (x_1,x_2),x_1^2+x_2^2<1\rbrace$.
\end{rem}

\paragraph{Computing with real symmetries}
The following proposition is of great use to work with real
symmetries.
\begin{prop}\label{liftanti}
Let $Q$ be an $\R$-plane, and $\sigma_Q$ be the symmetry about $Q$. There exists a matrix
$M_Q\in{\rm SU(2,1)}$ such that
\begin{equation}
M_Q\overline{M_Q}=1\mbox{ and } \sigma_Q(m)={\bf P}(M_Q.\bar\bm)\mbox{ for any $m\in\HdC$ with lift
$\bm$,}
\end{equation}
where {\bf P} : $\C^3\longrightarrow \C P^2$ is the projectivisation map.
\end{prop}
\begin{proof}
  In the special case where $Q=\HdR$, the identity matrix satisfy
  these conditions. In general, let ${\bf Q}$ be a lift to $\C^{2,1}$
  of $Q$, and $A$ be a matrix of SU(2,1) mapping $\R^3$ to ${\bf Q}$.
  The matrix $A{\bar A}^{-1}$ satisfies the above conditions.
\end{proof}

\begin{rem}\label{prodanti}
  Let $\sigma_1$ and $\sigma_2$ be real symmetries, with lifts $M_1$
  and $M_2$ given by Proposition \ref{liftanti}. The product
  $\sigma_1\sigma_2$ is a holomorphic isometry, and lifts to the
  matrix $M_1\overline M_2$. Similarly, if $h$ is a holomorphic
  isometry lifting to $H$, the conjugation $h\sigma_1 h^{-1}$ lifts to
  $H M_1\overline{H^{-1}}$.
\end{rem}

\n The isometry type of the product of two real symmetries is directly
related to the relative position of their mirrors. The following Lemma
is due to Falbel and Zocca in \cite{FZ} (see the next section for
information about the different isometry types).
\begin{lem}\label{FalZoc}
  Let $P_1$ and $P_2$ be two real planes, with respective symmetries
  $\sigma_1$ and $\sigma_2$.  Then
\begin{itemize}
\item The closures in $\HdC\cup\partial\HdC$ of $P_1$ and $P_2$ are
  disjoint if and only if the isometry $\sigma_1\sigma_2$ is
  loxodromic.
\item The intersection of the closures in $\HdC\cup\partial\HdC$ of
  $P_1$ and $P_2$ contains exactly one boundary point if and only if
  $\sigma_1\sigma_2$ is parabolic.
\item The intersection of the closures in $\HdC\cup\partial\HdC$ of
  $P_1$ and $P_2$ contains at least one point of $\HdC$ if and only if
  $\sigma_1\sigma_2$ is elliptic.
\end{itemize}
\end{lem}

\subsection{Classification of isometries.}
Let $A$ be a holomorphic isometry of $\HdC$. It is said to be 
\textit{elliptic} (resp. \textit{parabolic},
resp. \textit{loxodromic}) if it has a fixed point inside $\HdC$
(resp. a unique fixed point on $\partial\HdC$, resp. exactly two fixed
points on $\partial\HdC$). This exhausts all the possibilities.

\n Note that there is still a small
ambiguity among elliptic elements. An elliptic isometry will be
called a \textit{complex reflection} if one of its lifts to SU(2,1)
has two equal eigenvalues, else, it will be said to be \textit{regular
  elliptic}. 

\n As in the case of PSL(2,$\R$), there is an algebraic
criterion to determine the type of an isometry according to the trace
of one of its lifts to SU(2,1). An element of PU(2,1) admits three
lifts to SU(2,1) which are obtained one from another by multiplication
by a cube root of 1. Therefore its trace is well-defined up to
multiplication by a cube root of 1.

\begin{prop}\label{classauto}
  Let $f$ be the polynomial given by
  $f(z)=|z|^4-8\mbox{Re}(z^3)+18|z|^2-27$, and $h$ be a holomorphic
  isometry of $\HdC$.
\begin{itemize}
\item The isometry $h$ is loxodromic if and only if $f(\tr\, h)$ is positive.
\item The isometry $h$ is  regular elliptic if and only $f(\tr\, h)$ is negative.
\item If $f(h)=0$, then $h$ is either parabolic or a complex reflection.
\end{itemize}
\end{prop}
\begin{proof}
  Note that $f$ is invariant under multiplication of $z$ by a cube
  root of $1$.  The polynomial $f$ is the resultant of $\chi$
  and $\chi'$, where $\chi$ is the characteristic polynomial of a lift
  of $h$ to SU(2,1). See \cite{Go} (chapter 6) for details.
\end{proof}

\begin{rem}\label{trace3}
The function $f$ given in Proposition \ref{classauto} can be written in
 real coordinates as

$$f(x+iy)=y^{4}+y^{2}\left(x+6-3\sqrt{3}\right)\left(x+6+3\sqrt{2}\right)
+\left(x+1\right)\left(x-3\right)^{3},$$ \n with $x,y\in\R$. From this
writing of $f$, it follows that then $h$ is loxodromic whenever
$\mbox{Re}(\tr(h))>3$.
\end{rem}

\begin{rem}\label{spectrum}
  It is a direct consequence of the definition of SU(2,1) that the
  set of eigenvalues of a matrix $A\in$SU(2,1) is invariant under the
  transformation $z\longmapsto 1/\bar z$. We again refer the reader to
  \cite{Go} (chapter 6).
\end{rem}

\paragraph{Loxodromic isometries}
The following facts about loxodromic isometries will be needed later.
\begin{prop}\label{loxonormal}
Let $h\in$PU(2,1) be a loxodomic isometry. Then $h$ is conjugate in PU(2,1) to an isometry given by the
matrix in SU(2,1)
\begin{equation}
{\bf D}_\lambda=
\begin{bmatrix}
\lambda & 0 & 0\\
0 & \bar\lambda/\lambda & 0\\
0 & 0 & 1/\bar\lambda
\end{bmatrix}\mbox{ with }\lambda\in\C, |\lambda|\neq 1. \label{Dlambda}
\end{equation}
\end{prop}
\begin{proof}
  Since PU(2,1) acts doubly transitively on the boundary of $\HdC$,
  the isometry $h$ is conjugate to a loxodromic isometry fixing the
  two points $\infty$, and $[0,0]$. These two points lift respectively
  to the vectors $\begin{bmatrix} 1 & 0 & 0\end{bmatrix}^T$ and
  $\begin{bmatrix} 0 & 0 & 1\end{bmatrix}^T$.As a consequence of this,
  any lift of $h$ tu SU(2,1) (written in the canonical basis) must be
  diagonal, and as a consequence of Remark \ref{spectrum}, has the
  form given above.
\end{proof}

\n The family $\lbrace {\bf D}_t, t\in\R_{>0}\rbrace$ defines a 1-parameter
subgroup of PU(2,1) containing only matrices with real trace greater
or equal to 3.

\begin{defi}\label{1param}
  Let $R_\gamma$ the 1-parameter subgroup of PU(2,1) given by
  $g_\gamma^{-1}\lbrace {\bf D}_t, t>0\rbrace g_\gamma$, where $\gamma$
  is a geodesic in $\HdC$ and $g_\gamma$ is an isometry mapping the
  geodesic $\gamma$ to the geodesic connecting $\infty$ and $[0,0]$.
\end{defi}

The subgroup $R_\gamma$ does not depend on the choice of
$g_\gamma$. However, the parametrisation of $R_\gamma$ depends on this
choice. This small ambiguity will not be important in the rest of the
paper.

\begin{rem}
  Let us give another characterisation of $R_\gamma$. An isometry $A$
  belongs to $R_\gamma$ if and only if for any real plane $P$
  containing $\gamma$, $A$ preserves $P$ and the two connected
  components of $P\setminus\gamma$. Indeed, we may normalise the
  situation in such a way that $P$ is $\HdR$ and $\gamma$ is the
  geodesic connecting the two points with Heisenberg coordinates
  $[0,0]$ and $\infty$, in which case $R_\gamma=({\bf D}_t)_{t>0}$.
  The two connected components of $\HdR\setminus \gamma$ are $C^+$ and
  $C^-$, where, in horospherical coordinates, $C^{+}=\lbrace(x,0,u),
  x>0\mbox{ and }u>0\rbrace$ and $C^{-}=\lbrace(x,0,u), x<0\mbox{ and
  }u>0\rbrace$. Now, any isometry preserving $\HdR$, and fixing both
  $[0,0]$ and $\infty$ lifts to SU(2,1) as the diagonal matrix
  diag$(t,1,1/t)$ with real $t$. It is then a straightforward
  computation to check that such an isometry preserves the connected
  components $C^+$ and $C^-$ if and only if $t$ is positive, that is,
  if it belongs to $({\bf D}_t)_{t>0}$.
\end{rem}

\paragraph{Parabolic isometries.}
There are two main types of parabolic isometries: they can
be either \textit{unipotent} or \textit{screw-parabolic}.

\begin{enumerate}
\item Heisenberg translations are unipotent parabolics. They are
  conjugate to isometries associated to one of the above matrices $T_{[z,t]}$
  in SU(2,1), that correspond to those unipotent parabolics fixing
  $\infty\in \partial\HdC$.

$$T_{[z,t]}=\begin{bmatrix}1 & -\bar z\sqrt{2} & -|z|^2+it\\
0 & 1 & z\sqrt{2}\\
0 & 0 & 1
\end{bmatrix}\mbox{, with } z\in\C, t\in\R.$$
There are two PU(2,1)-conjugacy classes of Heisenberg translations.
The first one contains \textit{vertical} translations, which
correspond to $z=0$. In this case $T_{[0,t]}-Id$ is nilpotent of order
2. These isometries preserve a complex line. The parabolic elements
$T_{[0,t]}$ preserves the complex line $\HuC$.  The second conjugacy
class contains \textit{horizontal} translation. This is when $z\neq
0$, in which case $T_{[z,t]}-Id$ is nilpotent of order 3. These
isometries preserve a real plane, which is $\HdR$ in the case where
$z$ is real and $t$ vanishes.
\item Screw-parabolic isometries are
  conjugate to a product $h\circ r$, where $h$ is a vertical
  Heisenberg translation and $r$ a complex reflection about the
  invariant complex line of $h$.
\end{enumerate}

\section{Real ideal triangles.\label{RIT}}
\subsection{Ideal triangles}
An ideal triangle is an oriented triple of boundary points of $\HdC$.

\begin{defi}
  Let $(p_1,p_2,p_3)$ be an ideal triangle. The quantity
\begin{equation}\label{defCartan}\A\left(p_1,p_2,p_3\right)=\arg\left(-\la\bp_1,\bp_2\ra\la\bp_2,\bp_3\ra
\la\bp_3,\bp_1\ra\right)\end{equation}
does not depends on the choice of the lifts of the $p_i$'s, and is
called the Cartan invariant of the ideal triangle $(p_1,p_2,p_3)$.
\end{defi}

\n The Cartan invariant classifies the ideal triangles, as stated in
the following Proposition (see \cite{Go} chapter 7 for a proof).

\begin{prop}\label{Cartan}
The Cartan invariant enjoys the following properties
\begin{enumerate}
\item Two ideal triangles are identified by an element of PU(2,1)
  (resp. an antiholomorphic isometry) if and only if they have the
  same Cartan invariant (resp. opposite Cartan invariants).
\item An ideal triangle has Cartan invariant $\pm\pi/2$ (resp. $0$) if
  and only if it is contained in a complex line (resp. a real plane).
\end{enumerate}
\end{prop}

\begin{defi}
  We will call any ideal triangle contained in a real plane a
  \textit{real ideal triangle}.  
\end{defi}
\n Since the three points are contained
  in a real plane we will as well refer to the 2-simplex determined by
  three points on the boundary of a real plane as a real ideal
  triangle.

\begin{rem}\label{uniktriang}
  Up to isometry, there is a unique real ideal triangle, as shown by
  Proposition \ref{Cartan}. If $\tau$ and $\tau'$
  are two real ideal triangles, there are exactly two isometries
  mapping $\tau$ to $\tau'$, $\varphi$ and $\psi$. One of them
  (say $\varphi$) is holomorphic and the other antiholomorphic. More
  precisely, denoting by $\sigma$ the real symmetry about the real
  plane containing $\tau$, $\psi=\varphi\circ\sigma$.
\end{rem}

\subsection{The invariant of a pair of adjacent ideal real triangles.\label{invZ}}
We say that two real ideal triangles are \textit{adjacent} if they
have a common edge. All the pairs of real ideal triangles we consider
are \textbf{ordered}.

\begin{lem}\label{referenceRtriang}
  Let $\tau_1$ and $\tau_2$ be two adjacent ideal real triangles,
  sharing a geodesic $\gamma$ as an edge. There exists a unique
  complex number $z\in\C\setminus\lbrace -1,0\rbrace$ such that the
  ordered pair of real triangles $(\tau_1,\tau_2)$ is
  PU(2,1)-equivalent to the ordered pair of ideal real triangles
  $(\tau_0,\tau_z)$ given by the Heisenberg coordinates of its
  vertices by

\begin{equation}\label{refxalpha}
  \tau_0=\left(\infty,[-1,0],[0,0]\right)\mbox{ and }\tau_z=\left(\infty, [0,0],[z,0]\right)
\end{equation}
\end{lem}

\begin{proof}
  As shown by Proposition \ref{Cartan} and Remark \ref{uniktriang},
  there exists a unique holomorphic isometry $h$ mapping $\tau_1$ to
  $\tau_0$ and $\gamma$ to the geodesic connecting $\infty$ to
  $[0,0]$.  The isometry $h$ maps the triangle $\tau_2$ to an ideal
  triangle of which vertices are a priori given in Heisenberg
  coordinates by $\infty$, $[0,0]$ and $[z,t]$ with $z\in\C$ and
  $t\in\R$. Using relation (\ref{lifthoro}), where $u=0$ sinve we are
  on the boundary, we lift the latter three points to the three vectors

$$\bm_\infty=\begin{bmatrix} 1 \\ 0 \\ 0\end{bmatrix},\, \bm_{0,0}=\begin{bmatrix} 0 \\ 0 \\ 1\end{bmatrix}\mbox{ and }\bm_{z,t}=\begin{bmatrix} -|z|^2+it \\ z\sqrt{2} \\ 1\end{bmatrix}.$$

The triple product is given by

$$\la\bm_\infty,\bm_{0,0}\ra\la\bm_{0,0},\bm_{z,t}\ra\la\bm_{z,t},\bm_\infty\ra=-|z|^2-it.$$
Using (\ref{defCartan}), we obtain
$\A(\infty,[0,0],[z,t])=\tan\left(t/|z|^2\right)$. Hence the Cartan
invariant of $h(\tau_2)$ vanishes if and only if $t=0$.
\end{proof}

\begin{rem}\label{deg}
  Notice that the special cases where $z$ equals to $0$
  and $-1$ correspond respectively to the degenerate cases where one
  of the triangles has two collapsed vertices, and where the
  two triangles have the same set of vertices.
\end{rem}

\begin{defi}
  Let $(\tau_1,\tau_2)$ be a pair of adjacent real ideal
  triangles. We will call the complex number $z$ associated to it by
  Lemma \ref{referenceRtriang} the invariant of the pair
  $(\tau_1,\tau_2)$, and denote it by ${\tt
    Z}(\tau_1,\tau_2)$.
\end{defi}

\begin{rem}\label{uniktriang2}
  Let $(a,b,c)$ be a real ideal triangle, and $z$ be a complex number
  different from $0$ and $-1$. From Proposition
  \ref{referenceRtriang}, we see that there exists a unique point $d$
  in $\partial\HdC$ such that $(a,c,d)$ is a real ideal triangle, and
  $\tZ\left((a,b,c),(a,c,d)\right)=z$.
\end{rem}

\begin{rem}\label{simKR}
  It is possible to give another description of the invariant ${\tt
    Z}$ of a pair of ideal triangles.  Let $p_1$, $p_2$, $p_3$ and
  $p_4$ be four points in $\partial\HdC$, such that
  $\tau_1=(p_1,p_2,p_3)$ and $\tau_2=(p_3,p_4,p_1)$ are two real
  ideal triangles. Let $C_{13}$ be the (unique) complex line
  containing $p_1$ and $p_3$. Neither $p_2$ nor $p_4$ belong to
  $C_{13}$, since the two corresponding ideal triangles are real. The
  complex line $C_{13}$ lifts to $\C^3$ as a complex plane. Let
  $\bc_{13}$ be a vector in $\C^{2,1}$ Hermitian orthogonal to this
  complex plane. Then $C_{13}=\bf{P}\left(\bc_{13}^\perp\right)$. Let
  $\bp_i$ be a lift of $p_i$ for $i=1,2,3,4$.  The invariant ${\tt
    Z}(\tau_1,\tau_2)$ is given by

\begin{equation}\label{formuleZ}
{\tt Z}(\tau_1,\tau_2)=-
\dfrac{\la \bp_4,\bc_{13}\ra\la\bp_2,\bp_1\ra}{\la\bp_2,\bc_{13}\ra\la\bp_4,\bp_1\ra}\end{equation}

\n The above quantity does not depend on the various choices of lifts we
made. This definition is similar to the one of the complex cross-ratio
of Kor\'anyi and Reimann (see \cite{KR}).  To check that this formula is
valid, it is sufficient to check it on the special case $p_1=\infty$,
$p_2=[-1,0]$, $p_3=[0,0]$ and $p_4=[z,0]$. In this case, the choice
$\bc_{13}=\begin{bmatrix} 0 & 1 & 0\end{bmatrix}^T$ is convenient.
This invariant is similar to the one used by Falbel in \cite{F},
although the form (\ref{formuleZ}) is not used there. Note that 
(\ref{formuleZ}) shows that ${\tt Z}$ is preserved by holomorphic
isometries.
\end{rem}

\begin{lem}\label{antiholoangle}
  Let $(\tau_1,\tau_2)$ be a pair of adjacent real ideal
  triangles, with ${\tt Z}(\tau_1,\tau_2)=z$, and $f$ be an
  antiholomorphic isometry. Then ${\tt
    Z}(f(\tau_1),f(\tau_2))=\bar z$.
\end{lem}
\begin{proof}
  Let $\sigma$ be the symmetry about the real plane containing
  $\tau_1$. The isometry $f\circ\sigma$ is holomorphic, and therefore
  preserves the invariant of pairs of adjacent real ideal triangles.
  As a consequence, it is sufficient to show that ${\tt
    Z}(\sigma(\tau_1),\sigma(\tau_2))=\bar z$. We can normalise
  the situation to the reference pair $(\tau_0,\tau_{z})$ given by
  (\ref{refxalpha}). In this case, the real symmetry $\sigma$ is just
  the complex conjugation. It fixes the three points $\infty$,
  $[-1,0]$ and $[0,0]$, and maps the point $[z,0]$ to $[\bar z,0]$.
\end{proof}

\begin{prop}\label{Rtrigsym}
  Let $\tau_1=(a,b,c)$ and $\tau_2=(a,c,d)$ be two adjacent real ideal
  triangles. There exists a unique real symmetry $\sigma$ such that
  $\sigma(a)=c$ and $\sigma(b)=d$.
\end{prop}
\begin{proof}
  It is sufficient to prove that such a real symmetry exists and is
  unique for standard case where the two triangles are $\tau_0$ and
  $\tau_{z}$. More precisely, we have to show that for any $z\in\C$
  and, there exists a unique real symmetry $\sigma_{z}$ such that

\begin{equation}\sigma_{z}([-1,0])=[z,0]\mbox{ and }
\sigma_{z}(\infty)=[0,0].\label{symetrie}\end{equation}

\n If there existed two such symmetries, their product would be a
holomorphic isometry having four fixed points on $\partial\HdC$, not
belonging to the boundary of a complex line. Therefore this product
would be the identity. Therefore such a symmetry is unique if it exists.

\n Writing $z=xe^{i\alpha}$, the matrix 
\begin{equation}\label{reflection}M_{x,\alpha}=\begin{bmatrix}0 & 0 & x\\ 0 &
    e^{i\alpha} & 0\\ 1/x & 0 &0\end{bmatrix}\end{equation} is such
that $M_{x,\alpha}\overline{M_{x,\alpha}}=1$, and the real symmetry
associated to it satisfies to (\ref{symetrie}).
\end{proof}

\begin{defi}\label{symetriepaire}
  We call the involution provided by Proposition \ref{Rtrigsym} the
  symmetry of the pair $(\Delta_1,\Delta_2)$ and denote it by
  $\sigma_{\Delta_1,\Delta_2}$.
\end{defi}

\begin{rem}\label{orient}
  As a direct consequence of Lemma \ref{antiholoangle} and Proposition
  \ref{Rtrigsym}, we see that for any pair $(\tau_1,\tau_2)$ of
  adjacent real ideal triangles,

$${\tt Z}(\tau_1,\tau_2)=\overline{{\tt Z}(\tau_2,\tau_1)}.$$
\end{rem}

\begin{rem}\label{0pi}
  \n Let us consider the special case where $\tt Z$ is real.  Going
  back to Lemma \ref{referenceRtriang}, we see that if $z$ is real,
  the two triangles $\tau_0$ and $\tau_z$ are both contained in the
  standard real plane $\HdR$ since their vertices all have real
  coordinates. In horospherical coordinates, this real plane is
  nothing but a copy of the upper half-plane. Therefore $z$ is positive if
  and only if the two triangles are in the same connected component of
  of $\HdR\setminus\gamma$, where $\gamma$ is the common geodesic edge
  of the two triangles. Applying Lemma \ref{referenceRtriang}, we
    obtain in general that

\begin{itemize}
\item $\tt Z(\tau_1,\tau_2)$ is real if and only if $\tau_1$ and
  $\tau_2$ lie in a common real plane $P$.
\item $\tt Z(\tau_1,\tau_2)$ is positive (resp. negative) if and only
  if $\tau_1$ and $\tau_2$ lie in opposite (resp. the same) connected
  components of $P\setminus\gamma$.
\end{itemize}
\end{rem}

\begin{rem}\label{linkclassical}
  When four points $(p_i)_i$ belong to the boundary of the
  standard real plane $\HdR$, the invariant
  $\tZ((p_1,p_2,p_3),(p_1,p_3,p_4))$ is the classical cross-ratio in
  the upper-half plane, as can be checked from the embedding of the
  upper half plane in $\HdC$ given by (\ref{embedhalfsup}) in section
  \ref{total}.
\end{rem}

\begin{figure}
\begin{center}
\begin{pspicture}(-2,-2)(4,4)
\scalebox{2}{\includegraphics[width=5cm,height=2cm]{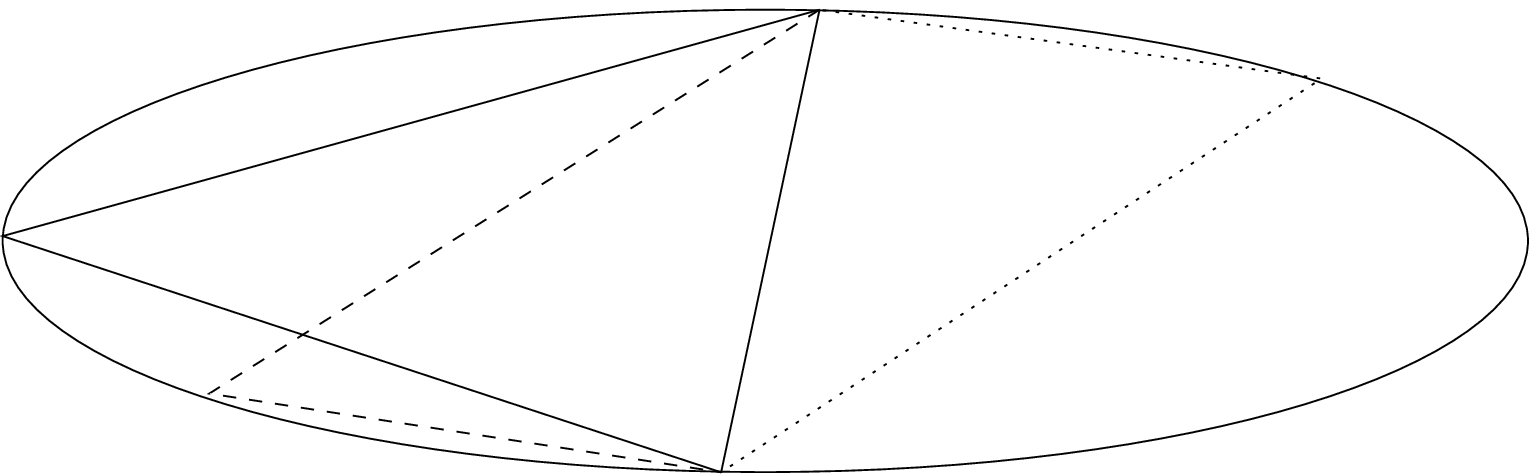}}
\rput(0,1){$P$}
\rput(-11,2){$\tau_1$}
\psline{->}(-7,.5)(-7,-0.5)
\rput(-7,-1){$\tau_2$ for $\tt Z(\tau_1,\tau_2)<0$}
\psline{->}(-9.5,2)(-10.5,2)
\psline{->}(-2,3)(-1,4)
\rput(-0.5,4.5){$\tau_2$ for $\tt Z(\tau_1,\tau_2)>0$}
\end{pspicture}
\caption{When $\tt Z(\tau_1,\tau_2)$ is real.}
\end{center}
\end{figure}

\n We will need the following in section \ref{defispinal}.

\begin{prop}\label{isomRtriangle}
  Let $\tau=(a,b,c)$ be an ideal real triangle, and $\gamma$ the
  geodesic connecting $a$ and $c$. Let $\tau_1$ and $\tau_2$ be
  two other real ideal triangles adjacent to $\tau$ along $\gamma$.
  Assume moreover that the invariants of the pairs $(\tau,\tau_1)$
  and $(\tau,\tau_2)$ satisfy

  \begin{equation}\label{hypo}\dfrac{{\tt Z}(\tau,\tau_1)}{{\tt
        Z}(\tau,\tau_2)}\in\R_{>0}.\end{equation}

\n Call $d_i$ the vertex of $\tau_i$ different from $a$ and $c$, and
$Q_i$ the mirror of the real symmetry $\sigma_i$ given by the
Proposition \ref{Rtrigsym}, such that $\sigma_i(a)=c$ and
$\sigma_i(b)=d_i$. Then there exists a unique element $g\in R_\gamma$
such that $g(Q_1)=Q_2$.
\end{prop}

\begin{proof}
We may normalise the situation so that
\begin{equation}
  a=\infty,\,b=[-1,0],\,c=[0,0],d_1=[z_1,0]\mbox{ and } d_2=[z_2,0],
\end{equation}
\n where $z_i={\tt Z(\tau,\tau_i)}$. In this normalised situation, the
two real symmetries associated to $d_1$ and $d_2$ are $\sigma_{z_1}$
and $\sigma_{z_2}$. As in \ref{reflection}, they correspond to the
matrices $M_{z_i}$ given by

\begin{equation}
M_{z_i}=\begin{bmatrix} 0 & 0 & |z_i|\\0 & z_i/|z_i| & 0 \\ 1/|z_i| & 0 & 0
\end{bmatrix}
\end{equation}

\n The one parameter subgroup $R_\gamma$ corresponds to $({\bf
  D}_t)_{t>0}$, with ${\bf D}_t$ as in Proposition \ref{loxonormal}).
Conjugating $M_{z_1}$ by ${\bf D}_t$ yields

\begin{equation}
{\bf D}_{t}M_{z_1}{\bf D}_{1/t}=\begin{bmatrix} 0 & 0 & |z_1|t^2\\
0 & \dfrac{z_1}{|z_1|} & 0\\
\dfrac{1}{|z|t^2} & 0 & 0\end{bmatrix}.
\end{equation}
Because of (\ref{hypo}), we have $z_1/|z_1|=z_2/|z_2|$, and the only possiblity is $t^2=z_2/z_1$, which leads to a unique value for $t$ since it is positive.
\end{proof}

\begin{rem}
  Keeping the notation and assumptions of Proposition
  \ref{isomRtriangle}, it is an easy exercice to check that in this
  case, the four points $a$, $b$, $d_1$ and $d_2$ belong to a common
  real plane, which is preserved by the isometry $g$. It is done by
  going back to the standard case of Lemma \ref{referenceRtriang},
  and looking at $\tau_0$, $\tau_{z_1}$ and $\tau_{z_2}$.
\end{rem}

\section{The bending theorem\label{sectionbending}}
\subsection{Notation, Definitions\label{decobend}}
We denote by $\Sigma=\Sigma_g\setminus\lbrace x_1,\cdots,x_n\rbrace$
an oriented surface of genus $g$ with $n$ deleted points such that
$2-2g-n<0$. We denote by $\pi_1$ its fundamental group, given by the
presentation
$$\pi_1=\la a_1,b_1,\dots a_g,b_g,c_1,\dots c_n\,\vert\,\prod_i[a_i,b_i]\prod_j c_j=1\ra,$$
where the $c_i$'s are homotopy classes of loops around the deleted
points.  The universal cover of $\Sigma$ is an open disk with a
$\pi_1$-invariant family of boundary points which may be thought of as
the lifts of the $x_i$'s. This family is called the \textit{Farey set}
of $\Sigma$, and we will denote it by $\mathcal{F}_\infty$. On way to
understand $\Fa$ is to endow $\Sigma$ with a finite area hyperbolic
structure. In this situation, $\Fa$ is the set of (parabolic) fixed
points of the $c_i$'s and their conjugates. If moreover the holonomy
of this hyperbolic structure is a subgroup of PSL($2$,$\mathbb{Z}$),
we obtain the classical Farey set $\mathbb{Q}\cup\infty$ in the
Poincar\'e upper half-plane.

Recall that an ideal triangulation of $\Sigma$ is a decomposition
$$\Sigma=\underset{\alpha}{\bigcup}\Delta_\alpha,$$
where each $\Delta_\alpha$ is homeomorphic to a triangle of which
vertices have been removed, and such that $\alpha\neq\beta\Rightarrow
\overset{\circ}{\Delta}_\alpha\cap\overset{\circ}{\Delta}_\beta=\emptyset$. It
is a classical fact using Euler characteristic that any ideal
triangulation of a surface of genus $g$ with $p$ deleted points has
$4g-4+2p$ triangles and $6g-6+3p$ edges.

\begin{defi}
  Let $T$ be a an ideal triangulation of $\Sigma$, and $\hat T$ be the
  associated triangulation of $\hat\Sigma$. We will call
  \textit{$\HdC$-realization bent along $T$}, or \textit{$T$-bent
    realization} of $\mathcal{F}_\infty(\Sigma)$ any pair
  $(\phi,\rho)$ such that
\begin{itemize}
\item $\rho$ is a representation $\pi_1(\Sigma)\longrightarrow$ Isom($\HdC$)
\item $\phi:\mathcal{F}_\infty\left(\Sigma\right)\longrightarrow
  \partial\HdC$ is a $(\pi_1(\Sigma),\rho)$-equivariant mapping.
\item for any face $\Delta$ of $\hat T$ with vertices $a$, $b$, and
  $c$, the three points $\phi(a)$, $\phi(b)$ and $\phi(c)$ are
  contained in the boundary of a real plane.
\end{itemize}
\end{defi}

\n The group Isom($\HdC$) acts on the set of $T$-bent realizations of
$\mathcal{F}_\infty$ by $g\cdot
\left(\phi,\rho\right)=\left(g\circ\phi,g\rho g^{-1}\right)$. We will
denote by $\mathcal{BR}_T$ the set of Isom($\HdC$)-classes of $T$-bent
realizations for this action.

\begin{figure}
\begin{center}
\begin{pspicture}(-2,-2)(4,4)
\scalebox{2}{\includegraphics[width=5cm,height=2cm]{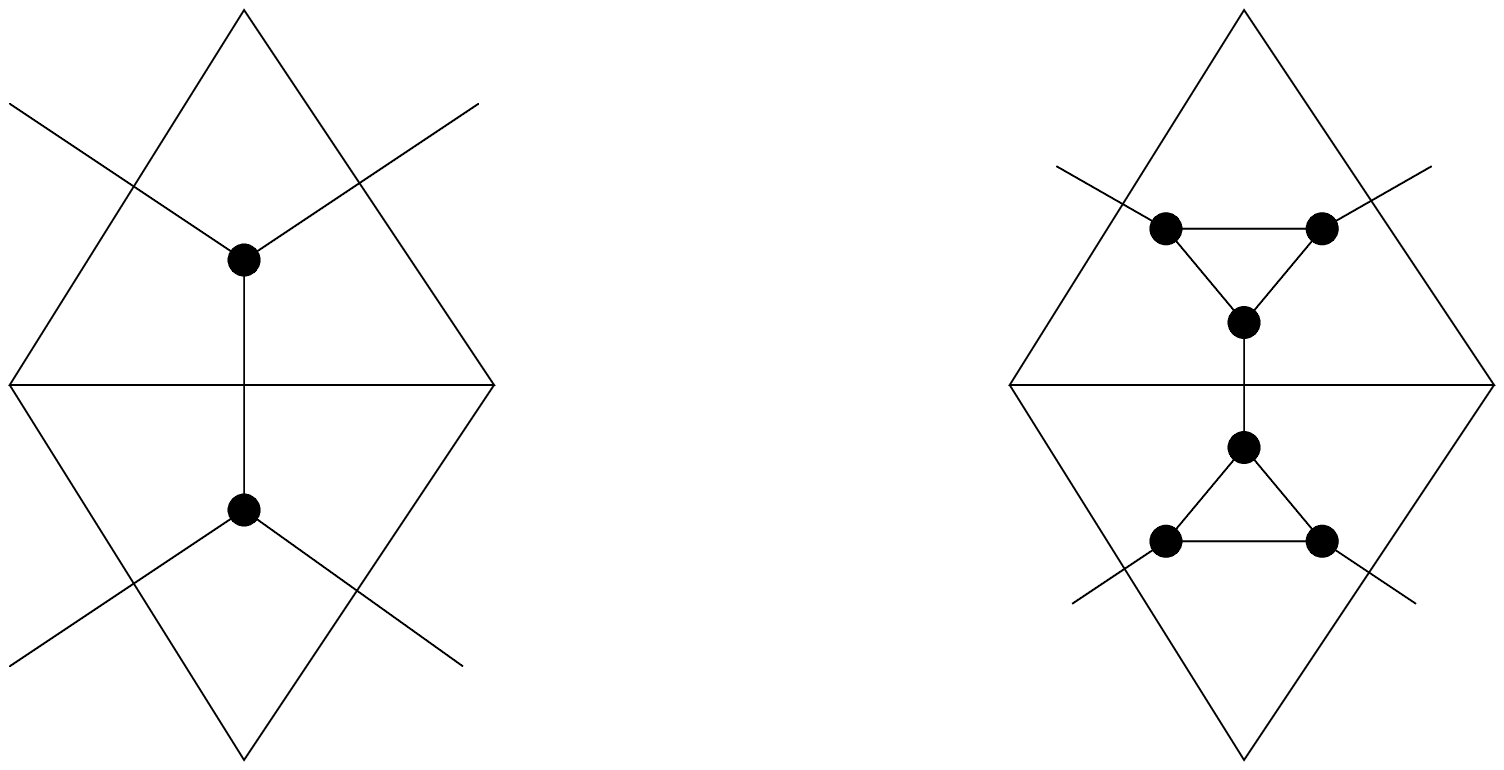}
\psarc{<-}(-2.5,0){1.5}{70}{110}}
\psline{<-}(-0.5,3.3)(-0.5,4.2)
\rput(0,4.5){\small{edge of type 1}}
\psline{->}(-1.3,1.5)(0.5,1.5)
\rput(1.7,1.5){\small{edge of type 2}}
\rput(0.8,2){$a_v=c_{v'}$}
\rput(-4,2){$a_{v'}=c_v$}
\rput(-1.7,4.3){$b_v$}
\rput(-1.7,-0.3){$b_{v'}$}
\rput(-2.1,2.4){$v$}
\rput(-2.1,1.8){$v'$}
\psarc{->}(-1.7,0.7){0.2}{100}{440}
\psarc{->}(-1.7,3.3){0.2}{100}{440}
\end{pspicture}
\caption{The triangulation, the modified dual graph and the labelling of vertices\label{modif}}
\end{center}
\end{figure}

\begin{defi}\label{dual}
  Let $T$ be an ideal triangulation of $\Sigma$. We will call
  \textit{modified dual graph} of T and denote by $\Gamma(T)$ the
  graph obtained from the dual graph of $T$ as follows (see figure
  \ref{modif}):
\begin{itemize}
\item The vertices of $\Gamma(T)$ are the combinations $1/3x+2/3y$,
  where $x$ and $y$ are adjacent vertices of the dual graph.
\item Two vertices $v$ and $v'$ of $\Gamma(T)$ are connected by an
  edge if and only if they fall in one of the following two cases.
\begin{itemize}
\item $v=1/3x+2/3y$ and $v'=2/3x+1/3y$ for some adjacent vertices $x$
  and $y$ of the dual graph. In this case the edge connecting $v$ and
  $v'$ is said to be of type 1.
\item $v=1/3x+2/3y$ and $v'=1/3z + 2/3 y$ where $yx$ and $yz$ are
  edges of the dual graph sharing an endpoint. In this case the edge
  $vv'$ is of type 2.
\end{itemize}
\end{itemize}
\end{defi}

We define similarly $\Gamma(\hat T)$, the modified dual graph of $\hat
T$, which is the lift of $\Gamma(T)$ to the universal cover of
$\Sigma$. We will refer to these two modified dual graphs as $\Gamma$
and $\hat\Gamma$ whenever it is clear from the context which
triangulation we are dealing with. Edges of type 1 and 2 of
$\hat\Gamma$ are defined similarly as for $\Gamma$. Note that an edge
of $\Gamma$ is of type 1 (resp. type 2) if and only if it intersects
an edge of T (resp. no edge of $T$). The orientation of $\Sigma$
induces an orientation of edges of type 2 of $\Gamma$ and
$\hat\Gamma$.

\begin{defi}\label{labelsomT}
Let $v$ be a vertex of $\hat\Gamma$ and $\Delta$ be the unique face
  of $\hat T$ containing $v$. The orientation of $\Sigma$
  induces an orientation of the edges of $\Delta$, and we will call
  $a_v$ the ending vertex of the edge of $\Delta$ closest to $v$. We
  will then call $b_v$ and $c_v$ the two other vertices of $\Delta$,
  in such a way that the triple $(a_v,b_v,c_v)$ is positively
  oriented.
\end{defi}

Since three vertices of $\hat\Gamma$ are contained in $\Delta$, there
are three possible labellings of the vertices of a given $\Delta$.

\begin{defi}
  A \textit{bending decoration} of an ideal triangulation $T$ is an
  application $\tD:e(T)\longrightarrow \C\setminus\lbrace -1,0\rbrace$
  defined on the set of unoriented edges of $T$.
\end{defi}

\n It follows from Remark \ref{deg} in section \ref{invZ}, that the
cases where the invariant ${\tt Z}(\tau_1,\tau_2)$ of a pair of real
ideal triangles equals $0$ or $-1$ correspond to degenerate pairs of
triangles: it equals $0$ if and only if $\tau_2$ has two identical
vertices and $-1$ if and only if the two triangles are equal. We do
not consider these degenerate cases.

\n We will often refer to the function $\arg(\tD)$ as the
\textit{angular part} of the bending decoration. There is an action of
$\mathbb{Z}/2\mathbb{Z}$ on the set of bending decorations of $T$
which is given by the complex conjugation: if $\tD$ is a bending
decoration of $T$, the decoration $\overline{\tD}$ is given by
$\overline{\tD}(e)=\overline{\tD (e)}$ for any edge $e$ of $T$.

\begin{defi}
  For any ideal triangulation $T$ of $\Sigma$, we denote by
  $\mathcal{BD}_T$ the set of bending decorations of $T$, and by
  $\mathcal{BD}^*_T$ the quotient of $\mathcal{BD}_T$ by the action of
  $\mathbb{Z}/2\mathbb{Z}$ given above.
\end{defi}

\n The set of bending decoration $\mathcal{BD}_T$ of $T$ is thus a
copy of $(\C\setminus\lbrace -1,0 \rbrace)^{|e(T)|}$. Prior to proving
Theorem \ref{realbend}, we introduce the following isometries of
$\HdC$.

\begin{defi}\label{elemHdC}
We will refer to the following isometries as \textit{elementary isometries}.

\begin{itemize}
\item For any $z\in\C\setminus\lbrace 0,1\rbrace$, we will call $\sigma_z$
  the real symmetry given by \ref{reflection}, where $z=xe^{i\alpha}$.
\item $\mathcal{E}$ is the isometry given by its lift to U(2,1), where
  \begin{equation}\mathcal{E}=\begin{bmatrix}-1 & \sqrt{2} & 1\\-\sqrt{2} & 1 & 0\\
      1 & 0 & 0\end{bmatrix}.\label{Ebis}\end{equation} \n We will
  identify $\mathcal{E}$ and its lift (notice that $-\mathcal{E}$ is
  in SU(2,1))
\end{itemize}

\end{defi}

\n The symmetry $\sigma_{z}$ acts on the complex hyperbolic space by
$\sigma_{z}(m)={\bf P}\left(M_{z}\bar\bm\right)$ (see Proposition
\ref{liftanti}). In Heisenberg coordinates, its action on the boundary
$\partial\HdC$ is given by

$$\sigma_{z}([w,t])=\left[\dfrac{z\bar w}{|w|^4+t^2}\left(|w|^2-it\right),
  \dfrac{t |z|^2}{|w|^4+t^2}\right].$$

\n From this we see that (and it follows from Proposition
\ref{Rtrigsym} as well)

$$\sigma_z(\infty)=[0,0], \sigma_z([-1,0])=[z,0].$$

\n The isometry $\mathcal{E}$ is elliptic of order 3 and permutes
cyclically the three points $\infty$, $\lbrack -1,0\rbrack$ and
$\lbrack 0,0\rbrack$.

\subsection{Bent $\HdC$-realizations : proof of Theorem \ref{realbend}.\label{proof2}}
We are now going to prove that there is a bijection between
$\mathcal{BD}^*_T$ and $\mathcal{BR}_T$.

\begin{proof}[Proof of Theorem \ref{realbend}]
  We will associate to any bending decoration in
  $\mathcal{BR}_T$ a unique pair of PU(2,1)-classes of $T$-bent
  realizations of $\Fa$, which represent the same
  Isom($\HdC$)--class and correspond to conjugate bending decorations.
  We will first associate to any vertex $v$
  of the modified dual graph a bent realization $(\phi_v,\rho_v)$ of
  $\Fa$ by using $v$ as a basepoint. We will see a
  posteriori that we obtain this way two PU(2,1)-classes of
  realization which correspond to the same
  Isom($\HdC$)-class.\\

  \n\textbf{Step 1: Definition of the mapping $\phi_v$.} We would like
  to interpret the complex numbers $\tD(e)$ as invariants of pairs of
  real ideal triangles, and use it in order to construct $\phi_v$
  reccursively. Each edge $e$ belongs to two faces of $\hat T$, say
  $\Delta_1$ and $\Delta_2$, and we have to chose whether we see
  $\tD(e)$ as $\tZ(\phi(\Delta_1),\phi(\Delta_2))$ or as
  $\tZ(\phi(\Delta_2),\phi(\Delta_1))=\overline{\tZ(\phi(\Delta_1),\phi(\Delta_2))}$. We do it using a bicoloring of $\hat T$.

  Let $v$ be a vertex of $\hat\Gamma$ and $\Delta_v$ be the face
  containing $v$. Label by $a_v$, $b_v$ and $c_v$ the vertices of the
  face of the triangulation $v$ belongs to, as in Definition
  \ref{labelsomT}.

\begin{itemize}
\item Give to the face $\Delta_v$ the colour white, and define
  $\phi_v(a_v)=\infty$, $\phi_v(b_v)=[-1,0]$ and $\phi_v(c_v)=[0,0]$.
\item Colour all the faces of $\hat T$ in black or white from
  the one containing $v$ by following the rule that two triangles
  sharing an edge have opposite colour.
\item Define the images of all the other points of
  $\mathcal{F}_\infty$ recursively according to the following
  principle: if an edge $e$ separates two faces $\Delta_w$ (white) and
  $\Delta_b$ (black), then the number $z$ associated to the edge $e$
  is interpreted as the invariant of the (ordered) pair
  $\tZ(\phi_v(\Delta_w),\phi_v(\Delta_b))$. If $\phi_v(\Delta_w)$ is
  already constructed, this defines $\phi_v(\Delta_b)$ unambiguously,
  shown by Remark \ref{uniktriang2}.
\end{itemize}

\n\textbf{Step 2: Definition of the representation $\rho_v$.} For any
$\gamma\in \pi_1$, we have to define an isometry $g_\gamma$ such that
$\phi_{v}(\gamma\cdot m)=g_\gamma\phi_{v}(m)$ for any $m$ in $\Fa$. In particular, such an isometry must map the reference triangle
$(\infty,[-1,0],[0,0])$ to the ideal real triangle $(\phi_v(\gamma
\cdot a_v),\phi_v(\gamma\cdot b_v),\phi_v(\gamma\cdot c_v))$. This can
be done in two ways, as shown by Remark \ref{uniktriang}, using either
a holomorphic or an antiholomorphic isometry. We define
$\rho_v(\gamma)$ according to the following rule (recall that $v$
belongs to a white triangle).
\begin{itemize}
\item If $\gamma\cdot v$ belongs to a white triangle, define
  $\rho_v(\gamma)$ to be the unique holomorphic isometry mapping
  $(\infty,[-1,0],[0,0])$ to $(\phi_v(\gamma \cdot
  a_v),\phi_v(\gamma\cdot b_v),\phi_v(\gamma\cdot c_v))$.

\item If $\gamma\cdot v$ belongs to a black triangle, choose the
  antiholomorphic one.
\end{itemize}

\n\textbf{Step 3: $\phi_v$ is $(\pi_1,\rho_v)$-equivariant.} If
$\rho_v(\gamma)$ is holomorphic, it preserves the invariant
$\texttt{Z}$. As a consequence of the definition of $\phi_v$ and
$\rho_v$, the identity $\phi_v(\gamma\cdot m)=\rho_v(\gamma)\phi_v(m)$
holds for any $m$, and for any $\gamma$ such that $\rho_v(\gamma)$ is
holomorphic. If $\rho_v(\gamma)$ is antiholomorphic, then it
transforms the invariants of real ideal triangle from $z$ to $\bar z$,
as seen in Lemma \ref{antiholoangle}. The equivariance property in
this case is a direct consequence of the choice made in the
construction of $\phi_v$ to interpret the decoration as
$\texttt{Z}(\phi_v(\Delta_w),\phi_v(\Delta_b))$.\\

\n\textbf{Step 4: Description of the class of the realization
  $(\phi_v,\rho_v)$.}  Let us compare first the classes of $T$-bent
realizations associated to two vertices $v$ and $v'$ of an edge
$e$ of $\hat \Gamma$.

\begin{enumerate}
\item Assume first that $v$ and $v'$ belong to different faces
  $\Delta$ and $\Delta'$ of the triangulation, that is, $e$ is of type
  1. These two faces have opposite colours and $e$ intersects an edge
  of $\hat T$, which is decorated by some complex number $z$. The
  vertices of $\Delta$ are $a_v$, $b_v$ and $c_v$ as in Definition
  \ref{labelsomT}. Call $d_v$ the vertex of $\Delta'$ which is not a
  vertex of $\Delta$. Then, according to their definitions, $\phi_v$
  and $\phi_{v'}$ satisfy to

$$\begin{matrix}
  \phi_v(a_v)=\infty &,& \phi_v(b_v)=[-1,0] &,& \phi_v(c_v)=[0,0] &
  \mbox{and} &
  \phi_v(d_v)=[z,0]\\
  \\
  \phi_{v'}(a_v)=[0,0] &,& \phi_{v'}(b_v)=[z,0] &,&
  \phi_{v'}(c_v)=\infty & \mbox{and} &
  \phi_{v'}(d_v)=[-1,0].\\
\end{matrix}$$ 

The antiholomorphic involution $\sigma_{z}$ (Definition
\ref{elemHdC}), is the unique isometry exchanging $\infty$ and $[0,0]$
on one hand, and $[-1,0]$ and $[z,0]$. Therefore we see
$\phi_{v'}=\sigma_{z}\circ\phi_v$, and
$\rho_{v'}=\sigma_{z}\rho_v\sigma_{z}$, that is
$(\phi_{v'},\rho_{v'})=\sigma_z\cdot(\phi_v,\rho_v)$. In this case,
the two realizations are in the same isometry class, but not the
same PU(2,1)-class.
\item By examining similarly what happens when $v$ and $v'$ are
  connected by an edge of type 2, that is, if they belong to a common
  face of $\hat T$, we see that
  $(\phi_v,\rho_v)=\mathcal{E}\cdot(\phi_{v'},\rho_{v'})$ if the
  orientation induced on $e$ by the orientation of $\Sigma$ is
  $v\rightarrow v'$, and
  $(\phi_v,\rho_v)=\mathcal{E}^{-1}\cdot(\phi_{v'},\rho_{v'})$ in the
  opposite case. The two realizations have the same holomorphic class
  in this case.
\end{enumerate}

\n If $v$ and $v'$ are arbitrary vertices of $\hat T$, belonging to
the triangles $\Delta_v$ and $\Delta_{v'}$ of $\hat T$, colour the
faces of $\hat T$ starting from $\Delta_v$.  The facts {\it 1} and
{\it 2} above imply that
\begin{itemize}
\item if $\Delta_v$ and $\Delta_{v'}$ have the same colour for this
  choice of coloring, then $(\phi_v,\rho_v)$ and
  $(\phi_{v'},\rho_{v'})$ correspond to the same PU(2,1)-class of
  $T$-bent realization,
\item if not, then $(\phi_v,\rho_v)$ and $(\phi_{v'},\rho_{v'})$
  correspond to the same Isom($\HdC$)-class, but have opposite
  PU(2,1)-classes.
\end{itemize}
\n Indeed, if $\Delta_v$ and $\Delta_{v'}$ have the same colour if and
only if any simplicial path connecting $v$ and $v'$ contains an even
number of edges of type 1. Since the PU(2,1)-class changes
every time an edge of type 1 is used, this shows the above assertion.\\

\n\textbf{Step 5: Passing from $\tD$ to $\overline{\tD}$.} We have so
far associated to $\tD$ a pair of PU(2,1)-classes of $T$-bent
realizations.  The choice of a starting vertex $v$ of $\hat\Gamma$
determines a coloring of the faces of $\hat T$. Call $r_w$ the class
corresponding to white triangles for this choice of coloring, and
$r_b$ the one corresponding to black triangles. If we keep the same
starting vertex $v$ but construct the classes associated to the
decoration $\overline{\tD}$, the new equivariant mapping
$\psi_v:\mathcal{F}_\infty\longrightarrow\partial\HdC$ is defined
recursively from

$$\begin{matrix}\psi_v(a_v)=\infty &,& \psi_v(b_v)=[0,0] &,& \psi_v(c_v)=[-1,0]&\mbox{ and }&  \psi_v(d_v)=[\bar z,0].\end{matrix}$$

\n As a consequence, we see that $\psi_v=\sigma\circ\phi_v$, where
$\sigma$ is the complex conjugation. The corresponding holonomy
representation are conjugate by $\sigma$. Therefore the change
$\tD\longrightarrow\overline{\tD}$ induces the permutation 
$(r_w,r_b)\longrightarrow(r_b,r_w)$.\\

\n\textbf{Step 6: The reverse operation: decorating a triangulation
  from a $T$-bent realization.}  Let $r=(\phi,\rho)$ be a $T$-bent
realization of $\Fa$ in $\partial\HdC$. Since $r$ is bent along $T$ we
obtain by definition a family of real ideal triangles by connecting
$\phi(m)$ and $\phi(n)$ each time $m$ and $n$ are connected by an edge
of $\hat T$. If $e$ is an edge of $\hat T$ belonging to two triangles
$\Delta$ and $\Delta'$. As before colorating the the faces of $\hat T$
gives a way to associate to $e$ a complex number $z$, which is
$\tZ(\Delta,\Delta')$ if $\Delta$ is white and $\Delta'$ is black, and
$\overline\tZ(\Delta,\Delta')$ in the other case. There is an order 2
ambiguity: if we start with a given real ideal triangle, and obtain
this way a decoration $\tD$, starting with an adjacent triangle will
produce the decoration $\overline{\tD}$.
\end{proof}

\subsection{Explicit computation of the representations\label{explicit}}
In this section we assume that $\Gamma$ is endowed with a decoration $\tD$.
\begin{defi}
  For any oriented edge $\nu$ of $\Gamma$, let $A_\nu$ be the isometry
  defined as follows (see Definition \ref{elemHdC}).
\begin{enumerate}
\item If $\nu$ is of type one and intersects an edge $e$ of $\hat T$,
  then $A_\nu$ is the real symmetry $\sigma_{\tD(e)}$.
\item If $\nu$ is of type two, then if it is positively oriented with
  respect to the orientation of $\Sigma$, $A_\nu=\mathcal{E}$, else
  $A_\nu=\mathcal{E}^{-1}$.
\end{enumerate}
\end{defi}

\begin{prop}\label{isomcheminHdC}
  Let $T$ be an ideal triangulation of $\Sigma$, with a bending
  decoration, $v$ and $v'$ be two vertices of $\Gamma$, and
  $p_{v,v'}=s_1\cdots s_k$ be a simplicial path connecting them. Call
  $r_v$ and $r_{v'}$ the $T$-bent realizations associated to $v$ and
  $v'$, and $B_{v,v'}$ be the isometry $A_{s_1}\cdots A_{s_k}$. Then
  $B_{v,v'}$ satisfies to
$$r_v=B_{v,v'}\cdot r_{v'}.$$
\end{prop}

\begin{proof}
  This is a direct reccursion using the second step of the proof of
  Theorem \ref{realbend}.
\end{proof}

\n We now compute the representation in terms of the bending decoration.

\begin{prop}\label{calculrepreHdC}
  Let $\gamma$ be a homotopy class of loop on $\Sigma$, and $v$ be a
  vertex of $\Gamma$. We may represent $\gamma$ as a simplicial path
  starting at $v$ consisting of a sequence $e_1\cdots e_k$ of oriented
  edges of $\Gamma$. Associate to $\gamma$ the isometry
  $B_{v,\gamma\cdot v}=A_{e_1}\cdots A_{e_n}$. Then
\begin{enumerate}
\item The isometry $B_{v,\gamma\cdot v}$ does not depend on the choice of the
  simplicial loop representing $\gamma$.
\item The mapping $\gamma\longmapsto B_{v,\gamma\cdot v}$ is equal to
  the representation $\rho_v$.
\end{enumerate} 
\end{prop}

\begin{proof}
\begin{enumerate}
\item
  It is a classical fact $\Sigma$, equipped with an ideal
  triangulation can be retracted onto the dual graph of $T$. Therefore
  any loop $l$ on $\Sigma$ can be homotoped to a sequence of edges of
  the dual graph. Once a basepoint is fixed, this loop corresponds to a
  sequence of edges of the dual graph of $\hat T$, the triangulation
  of $\hat\Sigma$ coming from $T$. But the dual graph of $\hat T$ is a
  tree and therefore the sequence of edges representing $l$ is unique
  (if we assume that two consecutive edges are distinct).  Passing
  from the dual graph to the modified dual graph, we lose this
  uniqueness property. Indeed, let $\Delta$ be a triangle of $\hat
  T$, crossed by this unique sequence of edges of the dual graph. Let
  $v_1$, $v_2$ and $v_3$ be the three vertices of $\hat \Gamma$
  belonging to $\Delta$. Then the original loop can be homotoped to
  the simplicial path $v_1\rightarrow v_2\rightarrow v_3$ or to
  $v_1\rightarrow v_3$ (see Figure \ref{homotop}). However the
  isometries associated to these two sequences of edges are
  $\mathcal{E}$ and $\mathcal{E}^{-2}$ or $\mathcal{E}^{-1}$ and
  $\mathcal{E}^{2}$, according to the orientation. Since $\mathcal{E}$
  has order three this does not change the contribution of this part
  of the path to $B_{v,\gamma\cdot v}$.
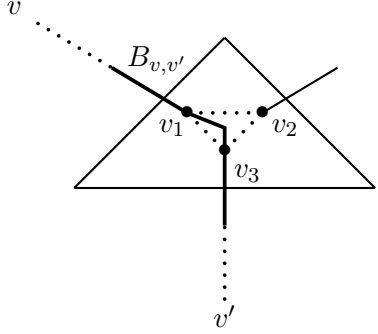
\begin{figure}
\begin{pspicture}(-2,-2)(2,2)
\psline(4,0)(6,2)
\psline(6,2)(8,0)
\psline(8,0)(4,0)
\psline[linestyle=dotted,linewidth=0.05](5.5,1)(6.5,1)
\psline[linestyle=dotted,linewidth=0.05](6.5,1)(6,0.5)
\psline[linestyle=dotted,linewidth=0.05](6,0.5)(5.5,1)
\rput(5.5,1){$\bullet$}
\rput(6.5,1){$\bullet$}
\rput(6,0.5){$\bullet$}
\rput(5.3,0.8){$v_1$}
\rput(6.8,0.8){$v_2$}
\rput(6.3,0.2){$v_3$}
\psline[linewidth=0.05](5.5,1)(4.5,1.6)
\psline[linestyle=dotted,linewidth=0.05](4.5,1.6)(3.5,2.2)
\rput(3.2,2.4){$v$}
\psline(6.5,1)(7.5,1.6)
\psline[linewidth=0.05](6,0.5)(6,-0.5)
\psline[linestyle=dotted,linewidth=0.05](6,-0.5)(6,-1.5)
\rput(6,-1.7){$v'$}
\rput(5.1,1.7){$B_{v,v'}$}
\psline[linewidth=0.05](5.5,1)(6,0.8)(6,0.5)
\end{pspicture}
\caption{Passing from the dual graph to the modified dual graph.\label{homotop}}
\end{figure}

 \item To prove the second assertion, we have to show that

\begin{enumerate}
\item $B_{v,\gamma\cdot v}$ maps the triple $(\phi_v(\gamma\cdot
  a_v),\phi_v(\gamma \cdot b_v), \phi_v(\gamma\cdot c_v)$ to the
  triple $(\infty,[-1,0],[0,0])$
\item $B_{v,\gamma\cdot v}$ is holomorphic if and only if $v$ and
$\gamma\cdot v$ lie in triangles having the same colour.
\end{enumerate}
 
We already know from Proposition \ref{isomcheminHdC} that
$\phi_{\gamma\cdot v}=A_{e_1}\cdots A_{e_n}\phi_{v}=B_{v,\gamma\cdot
  v}\phi_v$. As a consequence, the isometry $A_{v,\gamma\cdot v}$ maps
the triple $(\phi_v(\gamma\cdot a_v),\phi_v(\gamma \cdot b_v),
\phi_v(\gamma\cdot c_v)$ to the triple $(\phi_{\gamma\cdot
  v}(\gamma\cdot a_v),\phi_{\gamma\cdot v}(\gamma \cdot b_v),
\phi_{\gamma\cdot v} (\gamma\cdot c_v)$, which is by definition
$(\infty,[-1,0],[0,0])$. This shows the first part.

Now, the isometry $A_e$ attached to an edge $e$ is antiholomorphic if
and only if the edge $e$ is of type 1, that is, if $e$ passes from a
triangle to another. The isometry $B_{v,\gamma}$ is therefore
holomorphic if and only if the simplicial path corresponding to
$\gamma$ contains an even number of type 1 edges. Since the colour of
the triangle passes from black to white or vice versa at each
edge of type 1, we see that the isometry $B_{v,\gamma\cdot v}$ is
holomorphic if and only if the first and last triangles have the same
colour.
\end{enumerate} 
\end{proof}

\subsection{When is the representation in PU(2,1)?\label{when}}

\subsubsection{Representations in PU(2,1) and bipartite triangulations}
\begin{defi}\label{defibip}
  Let $T$ be an ideal triangulation of $\Sigma$, and $F$ be the set of
  faces of $T$. The triangulation $T$ is said to be bipartite if there
  exist two subsets of $F$, $F_1$ and $F_2$ such that
\begin{enumerate}
\item $F=F_1\bigcup F_2$
\item If a face $\Delta$ belongs to $F_i$, then its three neighbours
  belong to $F_{i+1}$, where the indices are taken modulo 2.
\end{enumerate}
\end{defi}

\begin{rem}
  An ideal triangulation is bipartite if and only if it is possible to
  colour its faces in two colours, black and white, in such a way
  that any white (resp. black) face has three black (resp. white)
  neighbours. For this reason, we will refer to black or white
  triangles.  Note that a triangulation is bipartite if and only if
  its dual graph is.
\end{rem}

\begin{rem}\label{relevbip}
  If $T$ is an ideal triangulation of $\Sigma$, then its lift $\hat T$
  to $\hat\Sigma$ is always bipartite. However, this bipartite structure of
  $\hat T$ projects onto a bipartite structure on $T$ if and only if
  it is $\pi_1$-invariant, that is, if and only if for any
  $\gamma\in\pi_1$ and any triangle $\Delta$ of $\hat T$, the two
  triangles $\Delta$ and $\gamma\cdot \Delta$ have the same colour.
\end{rem}

\begin{prop}\label{holobipar}
  Let $(T,\tD)$ be a decorated ideal triangulation of $\Sigma$, and
  let $\rho:\pi_1(\Sigma)\longrightarrow {\rm \textit{Isom}}(\HdC)$
  represent the Isom($\HdC$)-class of representation of $\pi_1$ in
  Isom($\HdC$) associated to $\tD$ by theorem \ref{realbend}. Then the
  following two statements are equivalent.
\begin{enumerate}
\item  The image of $\rho$ is contained in PU(2,1).
\item The triangulation $T$ is bipartite.
\end{enumerate}
\end{prop}
\begin{proof}

\begin{itemize}
\item We first prove that the bipartiteness is necessary. Pick a
  vertex $v$ of $\Gamma$ to be the basepoint. Let $\nu_1\cdots\nu_k$
  be a simplicial loop based at $v$ representing a homotopy class
  $\gamma\in\pi_1$.  Every $\nu_l$ of type 1 (resp. type 2)
  contributes to $\rho_v(\gamma)$ by an antiholomorphic (resp.
  holomorphic) isometry.  Hence $\rho_v(\gamma)$ is holomorphic if and
  only if $\nu_l$ is of type 1 for an even number of indices $l$. The
  number of colour changes is equal to the number of edges of type 1,
  and is even since $\gamma$ is a loop. Thus $\rho_v(\gamma)$ is
  holomorphic.

\item Assume now that $\rho(\gamma)$ is holomorphic for any
  $\gamma\in\pi_1$. Pick a homotopy class, and represent it by a
  simplicial loop $\gamma$ based at a vertex $v$ belonging to a face
  $\Delta_v$ of $T$. Attribute to $\Delta_v$ the colour white.  We can
  colour every triangles intersected by $\gamma$ by changing the
  colour every time an edge of type 1 is taken by $\gamma$. Since
  $\rho(\gamma)$ is holomorphic, the colour of $\Delta_v$ is
  well-defined (there are an even number of colour changes). We have to
  check now that if two simplicial loops $\gamma_1$ and $\gamma_2$
  based at $v$ intersect at a vertex $w\in\Delta_{w}$, then they
  define the same colour for $\Delta_{w}$.  Write these two loops

$$\gamma_1=\nu^1_1\cdots \nu^1_{k_1}\mbox{ and }
\gamma_2=\nu^2_1\cdots \nu^2_{k_2}.$$

Let $\gamma'_i$ one of the two subpathes of $\gamma_i$ connecting $v$
to $w$. Then $\gamma_{12}=\gamma'_1\gamma_{2}'^{-1}$ is a loop based
at $v$, and $\rho_v(\gamma_{12})$ is holomorphic. Therefore the number
of edges of type 1 in $\gamma_{12}$ is even. As a consequence, the
numbers of edges of type 1 in $\gamma'_1$ and $\gamma'_2$ have the
same parity and the colour of $\Delta_w$ is well-defined.
\end{itemize}
\end{proof}

\subsubsection{Existence of bipartite triangulations}
This section is devoted to the proof of the following proposition. 

\begin{prop}\label{existbipar}
  Let $\Sigma_{g,n}$ be a Riemann surface of genus $g$ with $n>0$
  deleted (or marked) points, such that $2-2g-n<0$.  Then
  $\Sigma_{g,n}$ admits a bipartite ideal triangulation.
\end{prop}

\begin{figure}
\begin{center}
  \begin{pspicture}(-2,-2)(2,2)
    \scalebox{0.5}{\includegraphics[width=15cm,height=5cm]{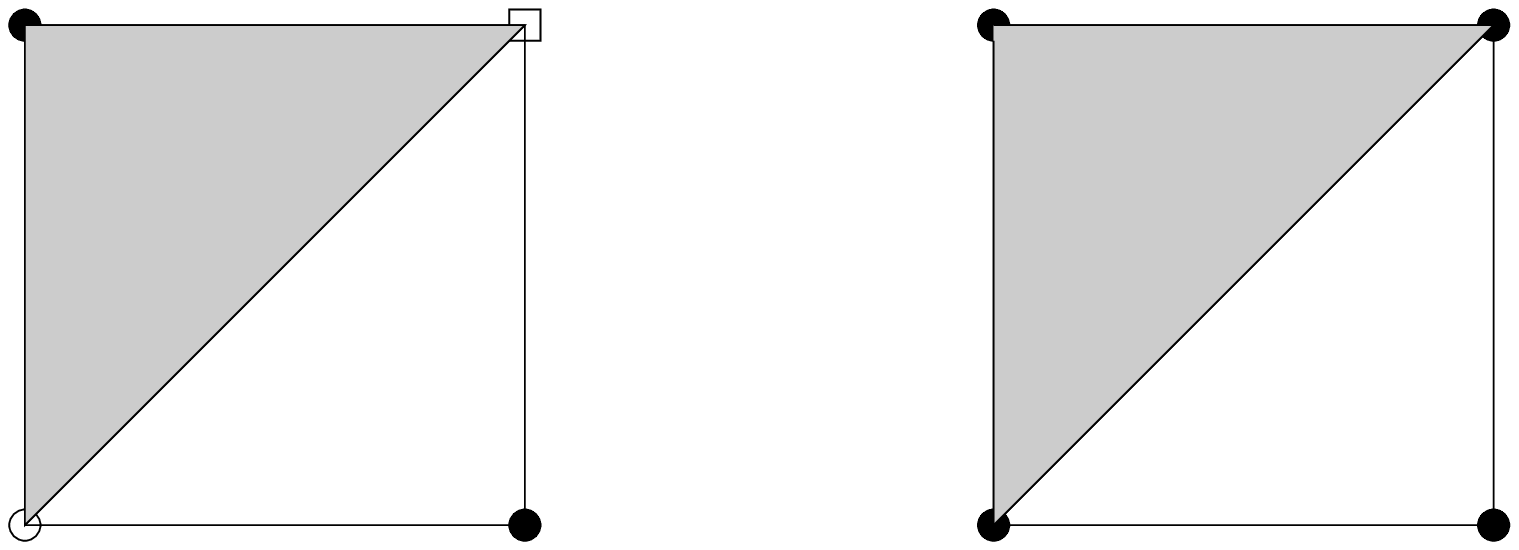}
      \psarc{->}(-15,0){1.5}{90}{360} \psarc{->}(-10,5){1.5}{-80}{180}
      \psline{->}(-2.5,0.3)(-2.5,4.7) \psline{->}(-4.9,2.5)(-0.3,2.5)}
    \rput(-6.5,-2){\small{The 3 -marked points sphere}}
    \rput(-0.5,-2){\small{The 1-marked point torus}}
\end{pspicture}
\caption{Bipartite ideal triangulations for surfaces of Euler characteritic -1\label{car-1}}
\end{center}
\end{figure}

\begin{figure}\vspace{2.2cm}
\begin{center}
  \begin{pspicture}(-2,-2)(2,2)
    \includegraphics[width=13cm,height=5cm]{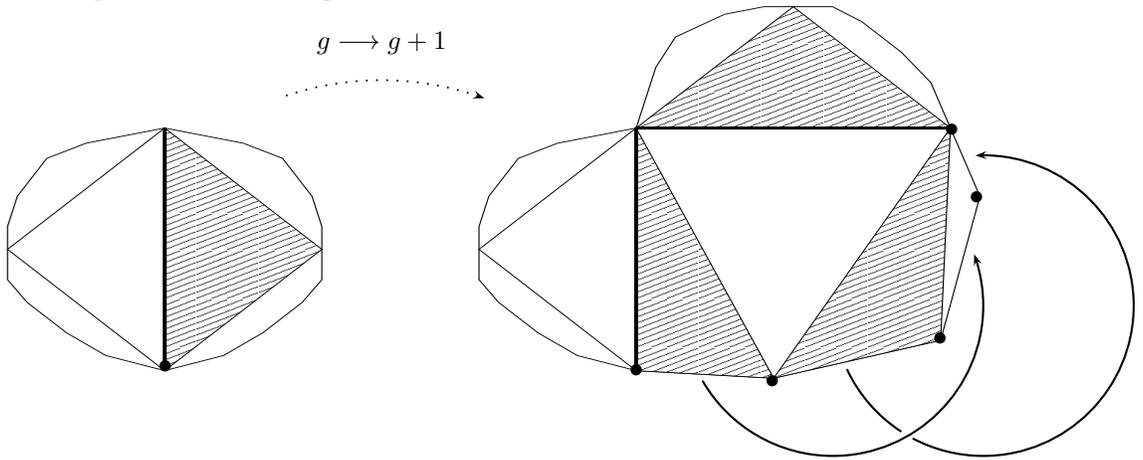}
    \put(-10.97,0.1){$\bullet$} \put(-4.71,0.05){$\bullet$}
    \put(-2.9,-0.1){$\bullet$} \put(-0.67,0.48){$\bullet$}
    \put(-0.19,2.36){$\bullet$} \put(-0.52,3.25){$\bullet$}
    \psarc[linestyle=dotted]{<-}(-8,0){4}{70}{110}
    \rput(-8,4.5){\small{$g\longrightarrow g+1$}}
    \psarc{->}(-2,1){2}{210}{380} \psarc(0,1){2}{205}{237}
    \psarc{->}(0,1){2}{242}{453}
  \end{pspicture}
  \caption{Increasing the genus\label{g+1}}
\end{center}
\end{figure}

\begin{figure}\vspace{2cm}
\begin{center}
  \begin{pspicture}(-2,-2)(2,2)
    \includegraphics[width=13cm,height=5cm]{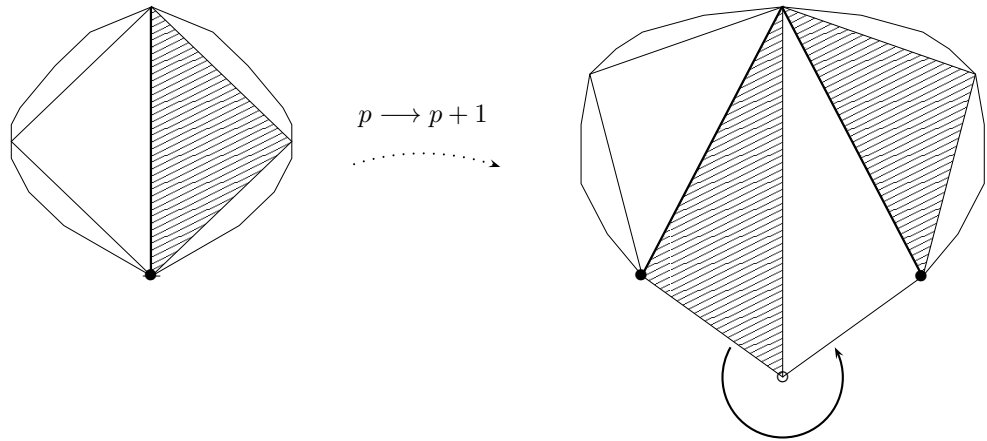}
    \put(-11.22,1.29){$\bullet$} \put(-4.7,1.29){$\bullet$}
    \put(-0.98,1.28){$\bullet$} \rput(-2.72,0.02){$\circ$}
    \psarc{->}(-2.72,0.02){0.8}{150}{390}
    \psarc[linestyle=dotted]{<-}(-7.5,0){3}{70}{110}
    \rput(-7.5,3.5){\small{$p\longrightarrow p+1$}}
  \end{pspicture}
  \caption{Increasing the number of marked points\label{p+1}}
\end{center}
\end{figure}

\begin{proof}
  We prove this proposition by induction, starting with the sphere
  with three marked points and the torus with one marked point.

  Both the 1-marked point torus and the 3-marked points sphere admit
  ideal triangulations consisting of two triangles, and the result
  is clear in these two cases (see Figure \ref{car-1}). We prove the
  result from these two cases by describing a recursion process
  increasing the genus of the surface by one or adding one puncture to
  the surface, and respecting the bipartiteness of the triangulation.
  We take the point of view that any triangulated surface is obtained
  from a triangulated polygon with identifications of the external
  edges.

  First, the bipartite triangulation of the surface corresponds to a
  bipartite triangulation of the polygon, compatible with the
  identification of external edges.  By this we mean that if two
  external edges are identified, then one of them should belong to a
  black triangle, and the other to a white one. We will denote
  respectively by $F$, $E$ and $V$ the sets of faces, edges and
  vertices of the triangulation.

\begin{itemize}
\item {\bf Increasing the genus} (see Figure \ref{g+1}). Pick an
  internal edge of the triangulated polygon, cut along it to open the
  polygon and insert four new triangles as on figure
  \ref{g+1}. Identify the new external edges created this way as
  indicated on figure \ref{g+1}. During this process, 4 new triangles
  were created, as well as 6 new edges and no new vertex. As a
  consequence, the Euler characteristic of the compactified surfaces
  changes from $\chi=|V|-|E|+|F|$ to
  $\chi'=|V|-(|E|+6)+(|F|+4)=\chi-2$. Since no new vertex was created,
  the genus has increased by 1. The bicoloring of the new polygon is
  compatible with the gluing. Therefore the corresponding
  triangulation of the surface is also bipartite.
\item {\bf Increasing the number of punctures} (see Figure \ref{p+1}).
  The method is the same, inserting this time two new triangles, as
  indicated on figure \ref{p+1}. This time the transformation changes
  $|V|$ to $|V|+1$, $|E|$ to $|E|+3$ and $|F|$ to $|F|+2$, and
  preserves $\chi$.  As a consequence, the genus of the surface does
  not change, and we have introduced a new deleted point on the
  surface.
\end{itemize}
\end{proof}

\subsection{Loops around holes\label{holes}}
\begin{figure}  
\begin{center}
\begin{pspicture}(-2,-2)(4,4)
\includegraphics[width=5cm,height=4cm]{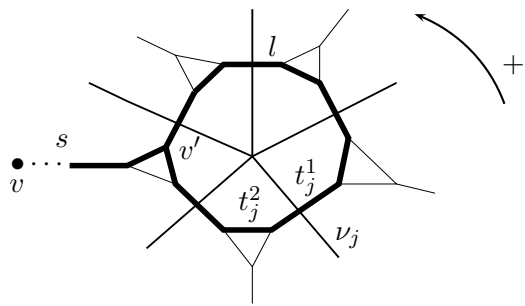}
\rput(-5.6,1.6){$v$}
\rput(-5.3,1.87){$\bullet\cdots$}
\rput(-3.3,2.1){$v'$}
\rput(-5,2.2){$s$}
\rput(-2.5,1.35){$t^2_j$}
\rput(-1.75,1.7){$t^1_j$}
\rput(-2.2,3.45){$l$}
\rput(-1.2,0.9){$\nu_j$}
\psarc{->}(-1,2){2}{20}{70}
\rput(1,3.2){+}
\end{pspicture}
\vspace{-2cm}
\caption{Loop around a vertex of the triangulation\label{autourducusp}}
\end{center}
\end{figure}

Let $T$ be a bipartite ideal triangulation of $\Sigma$ and $x$ be a
vertex of it. As seen in section \ref{proof2}, it is possible to
associate to any bending decoration of $T$ a class of $T$-bent
realisation $[\phi,\rho]$. We are going now to analyse the isometry
type of images of peripheral curves in terms of the bending decoration.

\begin{defi}\label{balanced}
  Let $T$ be an ideal triangulation of $\Sigma$, $\tD$ a bending
  decoration of $T$, $x$ one of the points deleted from $\Sigma$, and
  $\lbrace \nu_1,\cdots,\nu_k\rbrace$ be the set of edges of $T$
  having $x$ as an endpoint. We will say that $\tD$ is \textit{balanced at $x$} 
  whenever the following condition is satisfied

$$\prod_{i=1}^k|\tD(\nu_i)|=1.$$
We will say that a bending decoration is \textit{balanced} if it is
balanced at $x$ for every $x$.
\end{defi}

\n Notice that when $T$ is bipartite, the number of edges of $T$
having $x$ as a vertex is even. 

\begin{prop}\label{typecusp}
  Let $T$ be a bipartite ideal triangulation of $\Sigma$, $\tD$ be a
  bending decoration of $T$ and $r=(\phi,\rho)$ be the associated
  realisation of $\Fa$ in $\partial\HdC$. Let $c$ be a homotopy class
  of loop surrounding $x$ and no other puncture.  Then
\begin{enumerate}
\item The holomorphic isometry $\rho(c)$ is loxodromic if and only if
  $\tD$ is not balanced at $x$.
\item If $\tD$ is balanced at $x$, then the isometry $\rho(c)$ is
  either parabolic or a complex reflection.
\end{enumerate}
\end{prop}

\begin{proof}
  Pick a vertex $v$ of $\Gamma$ the modified dual graph of $T$ and
  represent the class $[\rho]$ by the representation $\rho_v$, as in
  Proposition \ref{calculrepreHdC}. Let $v'$ be a vertex of $\Gamma$
  belonging to one of the edges of $\Gamma$ intersecting one of the
  $\nu_j$'s (see Figure \ref{autourducusp}).  The homotopy class $c$
  is represented by a simplicial loop $sls^{-1}$, where $s$ in a
  simplicial path connecting $v$ to $v'$, and $l$ is a simplicial loop
  enclosing $p$, based at $v'$ such that $l=t^2_1t^1_1\cdots
  t^2_jt^1_j\cdots t^2_{k}t^1_k$, where $t_j^i$ is an edge of type $i$
  of $\hat \Gamma$ intersecting $\nu_j$ (see figure
  \ref{autourducusp}). Then, according to Proposition
  \ref{calculrepreHdC}, we see that $\rho_v(c)$ is conjugate to the
  product

\begin{equation}\sigma_{z_1}\circ \mathcal{E}^{\epsilon}\circ\sigma_{z_2}\circ
  \mathcal{E}^{\epsilon}\circ \cdots\circ\sigma_{z_{2k}}\circ
  \mathcal{E}^{\epsilon},\label{isomcusp}\end{equation}

\n where $z_j=\tD(\nu_j)$ and $\epsilon=1$ (resp. $-1$) when the
orientation of $c$ coincide with (resp. is opposite to) the one of the
surface.  The involution $\sigma_{z}$ being antiholomorphic, the
isometry (\ref{isomcusp}) products lifts to U(2,1) as the product of
matrices (see Remark \ref{prodanti}) 
\begin{equation}\label{prodcusp}
  M_{z_1}\mathcal{E}\overline{M_{z_2}}\mathcal{E}\cdots \overline{M_{z_{2k-1}}}\mathcal{E}M_{z_{2k}}\mathcal{E}=M_{z_1}\mathcal{E}M_{\bar z_2}\mathcal{E}\cdots
  M_{\bar z_{2k-1}}\mathcal{E}M_{z_{2k}}\mathcal{E}=\prod_{j=1}^{2p}M_{z_j^+}\mathcal{E},
\end{equation}
\n where $z_j^+$ is $z_j$ for odd $j$ and $\bar z_j$ for even $j$. For
any $z$, that the matrix $M_z\mathcal{E}$ is proportional to
the element of SU(2,1) given by
\begin{equation}\label{petitproduit}
  M_{z}\mathcal{E}\underset{\tiny{\mbox{SU(2,1)}}}{\sim}\begin{bmatrix}
    w & 0 & 0 \\
    -\sqrt{2}\bar w/w & \bar w/w & 0\\
    -1/\bar w & \sqrt{2}/\bar w & 1/\bar w\\
\end{bmatrix} \mbox{ where } w=\bar z^2/z.
\end{equation}

\n As a consequence, the product (\ref{prodcusp}) has diagonal
coefficients $\pi=\prod_{i=1}^{2p}w_i^+$, $\bar \pi/\pi$ and $1/\bar
\pi$ and is lower triangle.

The isometry $\rho(c_i)$ is therefore lowodromic if and only if the
product $\pi$ has modulus different from 1, that is, if
$\prod_{j=1}^{2k}|z_j|=\prod_{j=1}^{2k}|\tD(\nu_j)|\neq 1$ (notice
that $|w|=|z|$ in (\ref{petitproduit})). If $\pi$ has modulus 1, then
the isometry associated to the above matrix represents either a
parabolic isometry (if it is not semi-simple) or a complex reflection
(if it is semi-simple).

In the case where $\epsilon=-1$ is delt with in the same way,
with the only difference that $M_z\mathcal{E}^{-1}$ is upper
triangle instead of lower triangle. 
\end{proof}

\section{The discreteness theorem\label{sectiondiscret}}
\subsection{First part of the proof.\label{debut}}
The main goal of this section is to focus on those representations
associated to a special kind of bending decorations of the
triangulation, which we call \textit{regular}, and obtain our
discreteness results in this case.

\begin{defi}\label{decoregul}
  Let $T$ be a triangulation of $\Sigma$. We will say that a bending
  decoration $\tD$ of $T$ is \textit{regular} if there exists
  $\theta\in[-\pi,\pi[$ such that for all edges $e$ of $T$,
  $\arg(\tD(e))=\theta$.
\end{defi}

Recall that $x_1,\cdots,x_n$ are the points deleted from $\Sigma$, and
that $c_i$ denotes the class of peripheral loop surrounding $x_i$, in
the presentation of $\pi_1(\Sigma)$. Note that because of the
bipartiteness of $T$ the number of triangle having $x_i$ as a vertex
is even. Let us recall as well the statement of Theorem
\ref{theodiscret}.

\begin{theo*}[Theorem \ref{theodiscret}]
  Let $T$ be a bipartite ideal triangulation of $\Sigma$,
  $\theta\in]-\pi,\pi[$ be a real number and $\tD$ be a regular
  bending decoration of $T$ with angular part equal to $\theta$. Let
  $\rho$ be a representative of the (unique) Isom($\HdC$)-class of
  representations associated to $\tD$. Then
\begin{itemize}
\item For any index $i$, $\rho(c_i)$ is parabolic if and only if $\tD$
  is balanced at $x_i$.
\item The representation $\rho$ does not preserve any totally geodesic
  subspace of $\HdC$, unless $\theta=0$, in which case it is
  $\R$-Fuchsian.
\item As long as $\theta\in[-\pi/2,\pi/2]$, the representation $\rho$
  is discrete and faithful.
\end{itemize}
\end{theo*}

The first two parts of Theorem \ref{theodiscret} follows from what we
already know about bent representations. We will prove them now, and
postpone the proof of the last part of the result to section
\ref{fin}, after having introduced necessary material in section
\ref{defispinal}.

\begin{proof}[Proof of parts 1 and 2 of theorem \ref{theodiscret}]
\begin{enumerate}
\item To prove the first part of the theorem, let us go back to the
  proof of Proposition \ref{typecusp}. Consider $c$, one of the
  homotopy classes of loops around the holes, surroundind the deleted
  point $x$. Without loss of generality, we may assume that $c$ is
  positively oriented with respect to $\Sigma$.  Since $\tD$ is
  regular implies, $\rho(c)$ is conjugate to the isometry given by the
  following product of matrices, where $2k$ is the number of edges
  adjacent to $x$. 

\begin{equation}\label{produitreg}
M_{r_{1}e^{i\theta}}\mathcal{E}M_{r_2e^{-i\theta}}\mathcal{E}\cdots M_{
r_{2k}e^{-i\theta}}\mathcal{E}=\prod_{j=1}^{2k}M_{r_je^{(-1)^{j+1}i\theta}}\mathcal{E}, \mbox{ where }z_j=r_je^{i\theta}.
\end{equation}
Notice that the assumption about the orientation of $c$ implies that
$\epsilon=1$ in the proof of Proposition \ref{typecusp}. Let us be more
precise about this product. First, by a direct computation, we see
that

\begin{equation}M_{r_1e^{i\theta}}\mathcal{E}M_{r_2e^{-i\theta}}\mathcal{E}
=\begin{bmatrix}
r_1r_2 & 0 & 0\\
-\sqrt{2}(r_2e^{i\theta}+1) & 1 & 0\\
* & \dfrac{\sqrt{2}}{r_1r_2}(r_2e^{-i\theta}+1) & \dfrac{1}{r_1r_2}
\end{bmatrix}.\label{pour2}\end{equation}

Notice next that if we multiply (\ref{pour2}) on the right by a lower
triangle matrix $L$, the coefficients with indices $(2,1)$ and $(3,2)$
of the new matrix
$M_{r_1e^{i\theta}}\mathcal{E}M_{r_2e^{-i\theta}}\mathcal{E}L$ are
independant of the $*$ coefficient above.  Using this fact, it is a
straightforward reccursion to check that
the product (\ref{produitreg}) has the form

\begin{equation}\label{holocusp}
\begin{bmatrix}
  \prod_{j=1}^{2k}r_j & 0 & 0\\
  &&\\
  -A\sqrt{2} & 1 & 0\\
  &&\\
  * & -\bar A \sqrt{2}\prod_{j=1}^{2k}r_j^{-1}&
  \prod_{j=1}^{2k}r_j^{-1}
\end{bmatrix},\\
\end{equation}
\n where $z_j=r_ie^{i\theta}$ and
$$A=\underset{A_1}{\underbrace{1+\sum_{p=1}^{k-1}\prod_{j=2p+1}^{2k}r_j}} +e^{i\theta}
\underset{A_2}{\underbrace{\sum_{p=1}^{k}\prod_{j=2p}^{2k}r_j}}.$$

(See Remark \ref{k=3} below.)

\n The latter matrix corresponds to a loxodromic element if and only
if it has one eigenvalue of modulus greater than 1, that is if and
only if the product $\prod_{j=1}^{2k}r_j$ is different from 1. Thus
$\rho(c)$ is loxodromic if and only if $\tD$ is not balanced at $x$.

\n Assume now that $\prod_{j=1}^{2k}r_j=1$. Then the above matrix is
either the identity or a unipotent matrix in SU(2,1). If it were the
identity, $A$ would to be zero.

\begin{itemize}
\item If $e^{i\theta}$ is not real, $A$ is zero if and only if $A_1$
  and $A_2$ are. The positivity of the $r_i$'s implies that it is not
  the case.
\item If $e^{i\theta}$ is real, then $e^{i\theta}=1$ since we excluded
  the case where $\theta=\pi$. Again, the positivity of the $r_i$'s
  implies that $A$ is not zero in this case.
\end{itemize}

Therefore the product (\ref{produitreg}), and thus $\rho(c)$
 is unipotent if and only if $\tD$ is balanced at $x$.

\item For $\theta\in]-\pi,\pi[\setminus\{0\}$, we know by construction
  that any two adjacent ideal triangles are not contained in a common
  totally geodesic subspace of $\HdC$. Indeed, they cannot be in a
  complex line since each of them is real, and Remark \ref{0pi} implies
  that they are in a common real plane if and only $\theta$ is $0$ or
  $\pi$. But each of the vertices of the ideal triangles involved is
  the fixed point of a conjugate of one of the $\rho(c_i)$'s, all of
  which are non elliptic as checked above. The result follows then
  from Lemma \ref{projfix} below.
\end{enumerate}
\end{proof}

\begin{rem}\label{k=3}
  For the sake of lisibility, let us write down $A$ when $k=3$, that
  is if there are 6 edges adjacent to $x$. In this case:

$$A=1+r_5r_6+r_3r_4r_5r_6+e^{i\theta}\left(r_6+r_4r_5r_6+r_2r_3r_4r_5r_6\right).$$
\end{rem}

\begin{lem}\label{projfix}
  Let $\rho$ be a representation of $\pi_1(\Sigma)$ in PU(2,1)
  preserving a totally geodesic subspace $\mathcal V$ of $\HdC$, and
  such that none of the $c_i$'s is mapped to an elliptic
  isometry. Then all the fixed points of the $\rho(c_i)$'s belong to
  $\mathcal V$.
\end{lem}

\begin{proof}
  Call $p_\mathcal{V}$ the orthogonal projection onto $\mathcal V$,
  and let $c_i$ be such that $\rho(c_i)$ has a fixed point
  $m\in\partial\HdC\setminus\partial \mathcal V$. Since $\rho(c_i)$ is
  an isometry preserving $\mathcal V$, the two geodesic $(mp_{\mathcal
    V}(m))$ and $(m,\rho(c_i)(p_{\mathcal V}(m))$ are both orthogonal to
  $\mathcal V$. They are thus equal, and $\rho(c_i)$ fixes $p_{\mathcal
    V}(m)\in\HdC$, which is absurd since $\rho(c_i)$ is non-elliptic.
\end{proof}
\subsection{Spinal $\R$-surfaces.\label{defispinal}}
\n In order to prove the third part of the Theorem \ref{theodiscret},
we introduce in this section the main tool we will use.

\begin{defi}
  Let $P$ be an $\R$-plane, and $\gamma$ a geodesic contained in $P$.
  The \textit{spinal $\R$-surface} built on $\gamma$ with respect to
  $P$ is the hypersurface
$$S_{\gamma,P}=\Pi_P^{-1}\left(\gamma\right),$$
where $\Pi_P$ is the orthogonal projection onto $P$.
\end{defi}

\n Note that $\Pi_P$ is well-defined as the orthogonal projection onto
a totally geodesic subspace of a negatively curved Riemannian
manifold. It is a direct consequence of the definition that any two
spinal $\R$-surfaces are isometric, since PU(2,1) acts transitively on
the set of pairs $(\gamma,P)$, where $\gamma$ is a geodesic contained
in a real plane $P$. It is proved in \cite{PP1}, that if $P$ is a real
plane, $\sigma_P$ the symmetry about $P$ and $m$ a point of $\HdC$, a
lift to $\C^3$ of the projection of $m$ onto $P$ is given by

\begin{equation}\nonumber
   \dfrac{1}{|\bm|}\bm-\dfrac{\la\bm,\sigma_P(\bm)\ra}{|\la\bm,\sigma_P(\bm)\ra| |\sigma_P(\bm)|}\sigma_P(\bm),
\end{equation}

\n where $|\bm|=\sqrt{-\la\bm,\bm\ra}$. The above vector is
a representant of the midpoint of $m$ and $\sigma_P(m)$. In the
special case where $P$ is the standard real plane $\HdR$, 
$\sigma_P(\bm)=\bar\bm$ and $|\bm|=|\sigma_P(\bm)|$, and we obtain
as a lift of $\Pi_{\HdR}(m)$ to $\C^3$ the vector

\begin{equation}\label{liftprojHdR}
\bm-\dfrac{\la\bm,\bar\bm\ra}{|\la\bm,\bar\bm\ra|}\bar\bm.
\end{equation}

\n The latter expression of the projection extends to
$\HdC\cup\partial\HdC$.  We refer the reader to \cite{PP1} for more
information about this projection.

\begin{ex}
  Using the ball-model of $\HdC$, $\HdR$ is the real disc containing
  the points with real coordinates. Then the fibre of the orthogonal
  projection onto $\HdR$ over the point $(0,0)$ is the real plane
  $i\HdR=\lbrace (ix_1,ix_2),x_1^2+x_2^2<1\rbrace$. 
\end{ex}

\begin{rem}\label{bisector}
\begin{itemize}
\item In \cite{Most} (p. 185), Mostow defined \textit{spinal
    surfaces}, which are the inverse images of geodesics by the
  orthogonal projection onto a complex line instead of a real plane,
  or equivalently surfaces equidistant from two points in
  $\HdC$. Spinal surfaces are therefore foliated by complex
  lines. Note that if $\gamma$ is a geodesic, there exists a unique
  spinal surface containing it ($\gamma$ is referred to as its
  \textit{spine}).  In contrast, the set of spinal $\R$-surfaces
  containing a given geodesic $\gamma$ is parametrised by a circle
  $S^1$, since there is a circle of real planes containing $\gamma$.
\item Spinal $\R$-surfaces were already used in \cite{Wi2}, where they
  were called $\R$-balls.  They were then generalised to
  \textit{packs} by Parker and Platis in \cite{PP1}. In their
  terminology, spinal $\R$-surfaces correspond to \textit{flat packs}.
  The connection between packs and spinal $\R$-surfaces is given below
  by Lemma \ref{orbit}. See also a discussion in the survey
  \cite{PP2}.
\end{itemize}
\end{rem}
\begin{prop}\label{decriRsphere}
  The spinal $\R$-surface $S_{\gamma,P}$ is diffeomorphic to a ball of
  dimension 3, and is foliated by $\R$-planes. It separates $\HdC$ in
  two connected components which are exchanged by the symmetry about
  any of the leaves of the foliation.
\end{prop}
\begin{proof}
  The fibres of the orthogonal projection onto $P$ are $\R$-planes
  (see for instance \cite{PP1}).  Since $\R$-planes are discs, spinal
  $\R$-surfaces are diffeomorphic to $\R\times\HdR$, that is, a
  3-dimensional ball. A spinal $\R$-surface separates $\HdC$ in two
  connected components which are the inverse images of the two
  connected components of $P\setminus\gamma$ by the orthogonal
  projection onto $P$. Let $Q$ be a leaf of $S_{\gamma,P}$. We may 
  normalise so that in the ball model of $\HdC$, $P$ is $\HdR$,
  $\gamma$ connects the two points $(-1,0)$ and $(1,0)$,
  and $Q=i\HdR$. Then the symmetry about $Q$ acts on $\HdR$ by
  $(x_1,x_2)\longmapsto (-x_1,-x_2)$, and the two connected component
  are exchanged.
\end{proof}

\begin{prop}\label{interRsurfRplan}
  Let $\gamma\subset P$ be a geodesic contained in a real plane, and
  $P'$ be another real plane containing $\gamma$. Then one exactly of
  the following two possibilities occur.
\begin{enumerate}
\item The real plane $P'$ is contained in $S_{\gamma,P}$.
\item Each of the two connected components of $\HdC\setminus
  S_{\gamma,P}$ contains exactly one of the two connected components
  of $P'\setminus\gamma$.
\end{enumerate}
\end{prop}

\begin{proof}
  Let us use the ball model of $\HdC$.  Applying if necessary an
  element of PU(2,1), we may assume that
  $P=\HdR=\lbrace(x,y),x^2+y^2<1\rbrace$ and $\gamma=\lbrace
  (x,0),x\in]-1,1[\rbrace$. Denote by $\Pi$ the orthogonal projection
  onto $P$, and by $P^+$ (resp. $P^-$) the connected component of $P$
  containing points $(x,y)$ with $y>0$ (resp. $y<0$). Any real plane
  containing $\gamma$ is the image of $P$ under a rotation of angle
  $\alpha$ around $\gamma$, that is a transformation corresponding to

\begin{equation}
R_\alpha= 
\begin{bmatrix}
e^{-i\alpha/3} & 0 & 0 \\
0 & e^{2i\alpha_3} & 0\\
0 & 0 & e^{-i\alpha/3}
\end{bmatrix}\in \mbox{SU(2,1)},
\end{equation}
\n which acts in ball coordinates as
$(z_1,z_2)\longmapsto(z_1,e^{i\alpha}z_2)$. Note that $R_\alpha$ fixes
pointwise the complex line containing $\gamma$. We obtain this way a
family of  real planes $P_\alpha$ defined and parametrized by
$$P_\alpha=R_\alpha(\HdR)=\lbrace (x,e^{i\alpha}y),\,x^2+y^2<1\rbrace.$$

\n Note that $P_0=\HdR$. Since $P_{\alpha+\pi}=P_\alpha$, it is only
necessary to study the relative position of $P_\alpha$ and
$S_{\gamma,P}$ for $\alpha\in]0,\pi[$.  Let us pick a point
$m=(x,ye^{i\alpha})$ in $P_\alpha$. Following (\ref{liftprojHdR}), we
see that the projection of $m$ on $\HdR=P_0$ is given by the vector

\begin{figure}  
\begin{center}
\begin{pspicture}(-2,-2)(4,4)
\includegraphics[width=5cm,height=5cm]{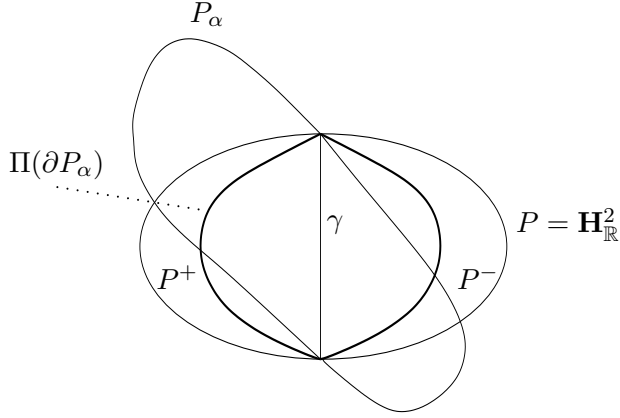}
\rput(0.8,2.5){$P=\HdR$}
\rput(-4,5.3){$P_\alpha$}
\rput(-2.3,2.5){$\gamma$}
\rput(-4.4,1.8){$P^+$}
\rput(-0.4,1.8){$P^-$}
\psline[linestyle=dotted](-4.1,2.7)(-6,3)
\rput(-6,3.3){$\Pi(\partial P_\alpha)$}
\end{pspicture}
\caption{Schematic picture for Proposition \ref{interRsurfRplan}.
\label{RsurfRplan}}
\end{center}
\end{figure}

\begin{equation}\label{projm}
\begin{bmatrix}
x \\ ye^{i\alpha}\\1
\end{bmatrix}
-
\dfrac{x^2+y^2e^{2i\alpha}-1}{|x^2+y^2e^{2i\alpha}-1|}\begin{bmatrix}
x \\ ye^{-i\alpha}\\1
\end{bmatrix}.
\end{equation}

\n Specialising (\ref{projm}) for $\alpha=\pi/2$, and using the fact
that $x^2-y^2-1\leq x^2+y^2-1<0$, we see that the second component of
the above vector vanishes. Thus any point $(x,iy)$ projects onto
$(x,0)$, which means that $P_{\pi/2}$ is contained in $S_{\gamma,P}$.  \\
\n We examine now the case where $\alpha\neq \pi/2$. Pick a point $m$
on $\partial P_{\alpha}$ distinct from $(\pm 1,0)$. The expression
(\ref{projm}) becomes (using $x^2+y^2=1$ and $\sin{\alpha}>0$):

\begin{equation}\label{projmbord}
\begin{bmatrix}
x \\ ye^{i\alpha}\\1
\end{bmatrix}
-
ie^{\alpha}\begin{bmatrix}
x \\ ye^{-i\alpha}\\1
\end{bmatrix}.
\end{equation}

\n which corresponds after projectivizing and rearranging to the point
with coordinates

$$\left(x,\dfrac{y\cos{\alpha}}{1+\sin\alpha}\right).$$
As a consequence, we see that a point $m=(x,ye^{i\alpha})\in \partial
P_\alpha$ with $y>0$ (resp. $y<0$) projects onto $P^+$ (resp. $P^-$)
if and only if $\alpha\in ]0,\pi/2[$, and that the situation is
opposite when $\alpha\in]\pi/2,\pi[$.  This proves the result for
connected components of the boundary of $P_\alpha$. \\

To conclude, let us assume that $\cos{\alpha}>0$ and that there is a
point $m=(x,ye^{i\alpha})$ in $P_\alpha$ with $y>0$ projecting onto
$P^-$. Then, considering the segment $\lbrace(x,te^{i\alpha}),
t\in[y,\sqrt{1-y^2}]\rbrace$ connecting $m$ to $\partial P_\alpha$,
we find a point with coordinates $(x,y'e^{\alpha})\in P_{\alpha}$
which projects to a point of $P$ with vanishing $y$ coordinate,
that is a point of $\gamma$. Applying if necessary a loxodromic
element in $R_\gamma$ (see Proposition \ref{loxonormal} and Definition
\ref{1param}), we may assume that the projection is actually the point
$(0,0)$. Since the fiber of $\Pi$ above $(0,0)$ is $i\HdR$, this
yields $\alpha=\pi/2$, which is absurd. The case where
$\cos{\alpha}<0$ is done in the same way. This proves the result
\end{proof} 

\n We give now another characterisation of spinal $\R$-surfaces.
Recall that if $\gamma$ is a geodesic, $R_\gamma$ is the 1-parameter
subgroup of PU(2,1) associated to $\gamma$. It contains the loxodromic
isometries of real trace greater than 3 preserving $\gamma$ (see
Definition \ref{1param}).

\begin{lem}\label{orbit}
  Let $Q$ be a real plane, and $\gamma$ be a geodesic of which
  endpoints we denote by $p$ and $q$. Assume that the real symmetry
  about $Q$ satisfies $\sigma_Q(p)=q$. Then the union $\cup_{g\in
    R_\gamma}g\cdot Q$ is a spinal $\R$-surface.  Conversely, any
  spinal $\R$-surface may be obtained in this way.
\end{lem}

\begin{proof}
  We may normalise the situation so that, using the ball model of
  $\HdC$, the points $p$ and $q$ have coordinates $p=(-1,0)$ and
  $q=(1,0)$, and $Q$ is the real plane $i\HdR$. The 1-parameter
  subgroup $R_\gamma$ preserves the real plane $\HdR$ and acts
  transitively on the geodesic connecting $p$ and $q$. Since $i\HdR$
  is the fibre of the orthogonal projection onto $\HdR$ above the
  point $(0,0)$ which belongs to $\gamma$, we see that $\cup_{g\in
    R_\gamma}g\cdot i\HdR$ is the spinal $\R$-surface built on
  $\gamma$ with respect to $P$.
\end{proof}

\n As said above, spinal surfaces enjoy two equivalent definitions,
either as inverse images of geodesics for the orthogonal projection
onto complex lines, or as surfaces equidistant from two given points
in $\HdC$. We have so far given an analogue of the first definition for
spinal $\R$-surfaces. The next proposition is more in the flavour of
the second one: it is possible to see spinal $\R$-surfaces as natural
obects separating two adjacent real ideal triangles, just as spinal
surfaces are naturally separating two distinct points. This version
of the definition will be of use to understand the geometric meaning
of the third part of Theorem \ref{theodiscret}.

\begin{prop}\label{separ}
  Let $\tau=(m_1,m_2,m_3)$ and $\tau'=(m_1,m_3,m_4)$ be two ideal
  real triangles and $\gamma$ be the geodesic connecting $m_1$ and
  $m_3$. Assume that the argument of ${\tt Z}(\tau,\tau')$ is not
  $\pi$. Then there exists a unique spinal $\R$-surface $S$ built on
  the geodesic $\gamma$ having the mirror of $\sigma_{\tau,\tau'}$
  as one of its leaves.
\end{prop}

\n Recall that $\sigma_{\tau,\tau'}$ is the symmetry of the pair
$(\tau,\tau')$ (see Definition \ref{symetriepaire}).

\begin{proof}
  Let $P$ be the mirror of $\sigma_{\tau,\tau'}$. Applying 
  Lemma \ref{orbit} to the real plane $P$ and the geodesic $\gamma$,
  we obtain a spinal $\R$-surface having the requested property. If
  there were another spinal $\R$-surface having the same property, the
  uniqueness part in Lemma \ref{Rtrigsym} would show that it would
  have $P$ as a leaf, and contain $\gamma$. Thus it would be equal to
  $S$ by Lemma \ref{orbit}.
\end{proof}
 
\begin{defi}\label{split}
Let $\tau$ and $\tau$ be two real ideal triangles sharing an edge and 
such that the argument of ${\tt Z}(\tau,\tau')$ is not $\pi$. 
We will call the spinal $\R$-surface given by Proposition \ref{separ} the 
\textit{splitting surface} of $\tau$ and $\tau'$ and denote it by 
Spl$(\tau,\tau')$.
\end{defi}

\begin{rem}\label{symsplit}
The definition of the splitting surface implies directly that
 Spl$(\tau_1,\tau_2)$=Spl$(\tau_2,\tau_1)$.
\end{rem}

\begin{prop}\label{splitplits}
  Let $\tau$ and $\tau'$ be two adjacent ideal triangles such that
  $\tZ(\tau,\tau')$ has argument different from $\pi$. Then $\tau$
  and $\tau'$ belong to opposite connected components of
  $\HdC\setminus\mbox{Spl}(\tau,\tau')$.
\end{prop}

\begin{proof}
  Since two spinal $\R$-surfaces are isometric, we may normalise the
  situation in such a way that the common geodesic of $\tau$ and
  $\tau'$ is in ball coordinates $\gamma=\lbrace
  (x,0),x\in]-1,1[\rbrace$, the splitting surface of $\tau$ and
  $\tau'$ is $S_{\gamma,\HdR}$, and the symmetry $\sigma_{\tau,\tau'}$ of the pair
  $(\tau,\tau')$ is the real symmetry about $i\HdR$, which is given in
  coordinates by

$$(z_1,z_2)\longmapsto(-\bar z_1,-\bar z_2).$$

We are in the same situation as in the proof of Proposition
\ref{interRsurfRplan}: $\tau$ is contained in one of the real planes
$P_\alpha$. Since $\tau'$ and $\tau$ are exchanged by
$\sigma_{\tau,\tau'}$, $\tau'$ is contained in the real plane
$\sigma_{\tau,\tau'}(P_\alpha)$, which is $P_{-\alpha}$. The result is
then a direct application of Proposition \ref{interRsurfRplan}.
\end{proof}

\begin{prop}\label{Rsurfangle}
  The splitting surface associated to a pair of
  adjacent real ideal triangles is determined by the argument of
  their ${\tt Z}$-invariant.
\end{prop}

\begin{proof}
  Let $\tau$ be a real ideal triangle, and $\gamma$ be one of its
  edges.  Consider $\tau_1$ and $\tau_2$ two real ideal triangles
  sharing the edge $\gamma$ with $\tau$ such that ${\tt
    Z}(\tau,\tau_j)=x_je^{i\alpha}$ for $j=1,2$. We have to show
  that the two spinal $\R$-surfaces Spl$(\tau,\tau_1)$ and
  Spl$(\tau,\tau_2)$ coincide.

  Call $Q_1$ and $Q_2$ the mirrors of the symmetries of the pairs
  $(\tau,\tau_1)$ and $(\tau,\tau_2)$. Proposition
  \ref{isomRtriangle} provides us a unique isometry $g$ belonging to
  the 1-parameter subgroup $G_\gamma$ which maps $Q_1$ to $Q_2$. In
  view of Lemma \ref{orbit}, the result is proved.
\end{proof}

\subsection{Proof of the third part of theorem \ref{theodiscret}\label{fin}}
We will prove now that a representation $\rho$ associated to a regular
bending decoration $\tD$ with anugular part $\theta\in[-\pi/2,\pi/2]$
is discrete and faithful. It is sufficient to prove that for these
values of $\theta$, the action of $\rho(\pi_1(\Sigma))$ acts properly
discontinuously on some $\rho(\pi_1(\Sigma))$-invariant subset of
$\HdC$. The following result is the crucial technical point.

\begin{theo}\label{crucial}
  Let $\tau$ be a real ideal triangle with vertices
  $(p_1,p_2,p_3)$. For $i=1,2,3$, let $\gamma_i$ be the geodesic
  $p_{i+1}p_{i+2}$ (indices taken mod. 3). Let $\tau_1$, $\tau_2$
  and $\tau_3$ be real ideal triangles, such that
\begin{itemize}
\item For $i=1,2,3$, $\tau$ and $\tau_i$ are adjacent, and share
  the geodesic $\gamma_i$ as an edge.
\item There exists $\theta\in [-\pi/2,\pi/2]$ such that
  $\arg\left({\tt Z}(\tau,\tau_i)\right)=\theta$ for $i=1,2,3$.
\end{itemize}
Then the three splitting surfaces $S_i={\rm Spl}(\tau,\tau_i)$ ($i=1,2,3$)
enjoy the following properties.

\begin{enumerate}
\item The intersection of $S_i$ and $S_{i+1}$ in $\HdC$ is empty.
\item The intersection of the closures of $S_i$ and $S_{i+1}$ in
  $\HdC\cup\partial\HdC$ is exactly $\lbrace p_{i+2}\rbrace$.
\end{enumerate}
\end{theo}

\n We postpone the proof of Theorem \ref{crucial}, and first finish
the proof of Theorem \ref{theodiscret}.

\begin{figure}  
\begin{center}
\begin{pspicture}(-2,-2)(4,4)
\includegraphics[width=5cm,height=5cm]{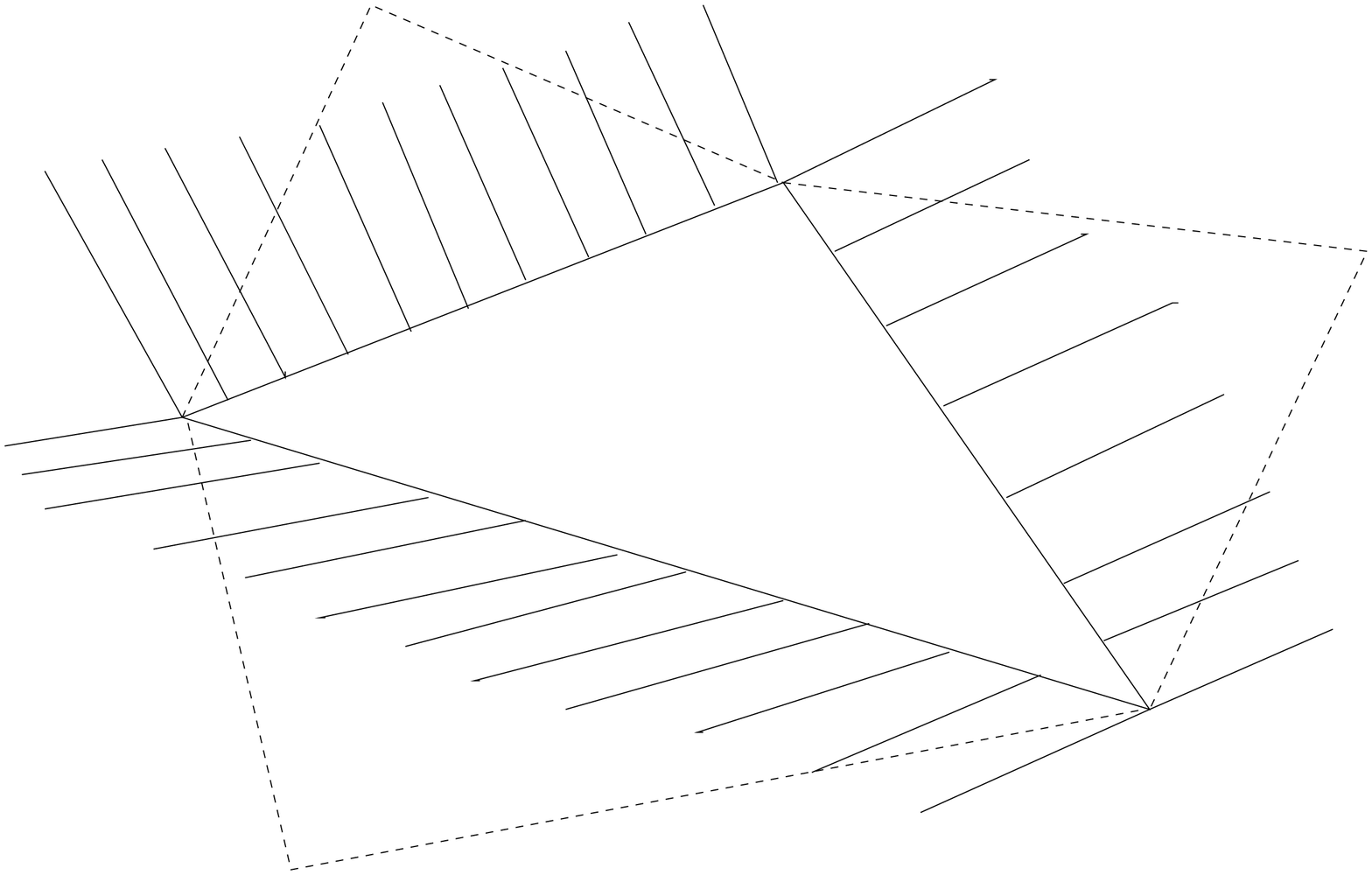}
\rput(-2.5,2.5){$\tau$}
\rput(0.1,3){$\tau_1$}
\rput(-4.4,0.2){$\tau_3$}
\rput(-4.1,5){$\tau_2$}
\rput(-3.5,4){\textbf{${S_2}$}}
\rput(-1,2.5){\textbf{${S_1}$}}
\rput(-2.7,1.3){\textbf{${S_3}$}}
\rput(-0.8,1){$\bullet$}
\rput(-4.3,2.65){$\bullet$}
\rput(-2.1,4){$\bullet$}
\rput(-0.6,0.7){$p_{2}$}
\rput(-2.2,3.6){$p_3$}
\rput(-3.9,2.65){$p_1$}
\end{pspicture}
\caption{Schematic picture for Theorem \ref{crucial}.\label{piccrucial}}
\end{center}
\end{figure}

\begin{proof}[Proof of part 3 of Theorem \ref{theodiscret}]
  Consider a given bipartite ideal triangulation $T$ of $\Sigma$
  and a regular bending decoration $\tD$ of $T$ with bending angle
  $\alpha$.  From $\tD$, we can construct as in the
  proof of Theorem \ref{realbend} a family of real ideal triangles in
  $\HdC$, and a $T$-bent realisation $(\phi,\rho)$ of $\Fa$ such that
  any two neighbouring triangles $(\Delta,\Delta')$ in $\hat T$ are
  mapped by $\phi$ to two real ideal triangle $(\tau,\tau')$ with a
  common edge and such that $\arg(\tZ(\tau,\tau'))=\pm\alpha$. We
  associate to each such pair of ideal triangle its splitting surface
  $Spl(\tau,\tau')$. Proposition \ref{splitplits} implies that $\tau$
  and $\tau'$ lie in opposite connected components of $\HdC \setminus
  Spl(\tau,\tau')$. Now, let $\tau$ be a triangle in the family, and
  $(\tau_i)_{i=1,2,3}$ be its three neighbours. Because of the
  regularity of the bending decoration, we have

\begin{equation}
  \arg{\left(\tZ(\tau,\tau_1)\right)}=\arg{\left(\tZ(\tau,\tau_2)\right)}
 =\arg{\left(\tZ(\tau,\tau_3)\right)}=\pm\alpha.
\end{equation}

Since $\alpha\in[-\pi/2,\pi/2]$, we can apply Theorem \ref{crucial},
and conclude that the three splitting surfaces $Spl(\tau,\tau_i)$ are
disjoint. Therefore each of the triangles $\tau$ obtained from the
bending decoration belong to a prism $\mathfrak{p}_\tau$, which is the
connected component of $\HdC\setminus \left(Spl(\tau,\tau_1)\cup
  Spl(\tau,\tau_2)\cup Spl(\tau,\tau_3)\right)$ whose boundary is made
of three spinal $\R$-surfaces. Applying reccursively Proposition
\ref{splitplits} shows that any two such prisms are either equal or
disjoint. As a consequence, we see that $\rho(\pi_1(\Sigma))$ acts on
the union of all the prims in such a way that

$$\rho(\gamma)\cdot \mathfrak{p}_\tau=\mathfrak{p}_{\rho(\gamma)(\tau)},\mbox{for any } \gamma \in \pi_1(\Sigma).$$
Therefore the action is discontinuous on $\cup_{\Delta\in\hat
  T}\mathfrak{p}_\Delta$, and $\rho$ is discrete.
\end{proof}

\n We prove now Theorem \ref{crucial}.

\begin{proof}[Proof of Theorem \ref{crucial}]$\mbox{(See figure \ref{piccrucial}})$.\\

  \n\textbf{First step: reduction to a normalised case.}\\
 By applying if necessary an isometry, we may assume that $\tau$ is
 the reference real ideal triangle given by $p_1=\infty$, $p_2=[-1,0]$
 and $p_3=[0,0]$. The isometry $\mathcal{E}$ given in by (\ref{Ebis})
 in Definition \ref{elemHdC} cyclically permutes the three latter
 points, and preserves the invariant $\tZ$ of pairs of real ideal
 triangles since it is holomorphic.  The bending decoration being
 regular, the invariants ${\tt Z}(\tau,\tau_{i})$ and ${\tt
   Z}(\tau,\tau_{j})$ have the same argument.  Therefore
 $\mathcal{E}$ maps $\tau_i$ to an ideal $\R$-triangle $\tau'_{i+1}$
 (indices taken mod. $3$) such that ${\tt Z}(\tau,\tau_{i+1})$ and
 ${\tt Z}(\tau,\tau'_{i+1})$ have the same argument. As a consequence
 of Proposition \ref{Rsurfangle}, it maps the splitting surface $S_i$
 to $S_{i+1}$, that is, it permutes the three splitting surfaces
 cyclically. Hence it is enough to prove that the two surfaces $S_1$
 and $S_2$ satisfy {\it 1} and
 {\it 2}.\\

  \n\textbf{Second step : parametrisation of the symmetries about the leaves of $S_2$ and $S_3$.}\\
  Let us use the following lifts for the $p_i$'s:
\begin{equation}
\bp_1=\begin{bmatrix}1 \\ 0\\ 0\end{bmatrix},
\bp_2=\begin{bmatrix}-1 \\ -\sqrt{2}\\ 1\end{bmatrix}\mbox{ and }
\bp_3=\begin{bmatrix}0 \\ 0\\ 1\end{bmatrix}.
\end{equation}

\n We first use Lemma \ref{orbit} to describe the leaves of $S_2$. Let
$q_2$ be the third point of $\tau_2$. According to Proposition
\ref{Rsurfangle}, we may assume that $q_2$ is any point such that
${\tt Z}(\tau,\tau_2)$ has the form $xe^{i\theta}$ with $x>0$. We
make the choice $q_2=[e^{i\theta},0]$.  The unique symmetry about a
real plane swapping $p_1$ and $p_3$, and $p_2$ and $q_2$ is given by
$\sigma_2(m)={\bf P}\left(M_2\bm\right)$, where $M_2$ is the matrix

$$M_2^\theta=\begin{bmatrix}
0 & 0 & 1\\
0 & e^{i\theta} & 0\\
1 & 0 & 0
\end{bmatrix}.$$

\n The 1-parameter subgroup $R_{\gamma_2}$ associated to the geodesic
connecting $p_1$ and $p_3$ is parametrised by the matrices

\begin{equation}
{\bf D}_{r_2}=\begin{bmatrix}
r_2 & 0 & 0\\
0 & 1 & 0\\
0 & 0 & 1/r_2
\end{bmatrix}\mbox{ with }r_2>0.
\end{equation}
\n We obtain thus the general form $M^\theta_{2,r_2}$ of a lift of the
symmetry about a leaf of $S_1$ by conjugating a lift of the involution
associated to $M_2^\theta$ by ${\bf D}_{r_2}$. Since $M_2^\theta$
stands for a antiholomorphic isometry, this yields (see Remark
\ref{prodanti})

\begin{eqnarray}
  M^\theta_{2,r_2} & = &{\bf D}_{r_2}M_{2}^\theta\overline{{\bf D}_{r_2}^{-1}}\nonumber\\
  & = & {\bf D}_{r_2}M_{2}^\theta{\bf D}_{1/r_2}\mbox{ (${\bf D}_{1/r_2}$ has real coefficients)}\nonumber\\
  & = &\begin{bmatrix}0 & 0 & r_2^2\\0 & e^{i\theta} & 0\\ 1/r_2^2 & 0 & 0\end{bmatrix}.
\end{eqnarray}

\n The general form $M^\theta_{3,r_3}$ of a lift of the symmetry about
a leaf of $S_3$ is obtained by conjugating the matrix
$M^\theta_{3,r_3}$ by the order three elliptic element $\mathcal{E}$:

\begin{eqnarray}
  M^\theta_{3,r_3} & = &\mathcal{E}{\bf D}_{r_3}M_{2}^\theta\overline{{\bf D}_{r_3}^{-1}}\overline{E^{-1}}\nonumber\\
  & = & \mathcal{E}{\bf D}_{r_3}M_{1}^\theta{\bf D}_{1/r_3}\mathcal{E}^{-1}\nonumber\\
  & = & \begin{bmatrix}
    -r_3^2& \sqrt{2}\left(e^{i\theta}+r_3^2\right)& \dfrac{1+2e^{i\theta}r_3^2+r_3^4}{r_3^2}\\
    &&\\
    \sqrt{2}r_3^2 & e^{i\theta}+2r_3^2& \sqrt{2}\left(e^{i\theta}+r_3^2\right)\\
    &&\\
    r_3^2 & -\sqrt{2}r_3^2 &  -r_3^2
\end{bmatrix}.
\end{eqnarray}
\n\textbf{Third step: proof of the disjunction}\\
Note first that the closures of $S_2$ and $S_3$ in
$\HdC\cup\partial\HdC$ both contain the point $p_1$ as a common end of
the geodesics $\gamma_2$ and $\gamma_3$. Therefore their intersection
should at least contain this point. Now, the result will be proved if
we show that the closure of any leaf of $S_2$ is disjoint from the
closure or any leaf of $S_3$.  We do this by showing that the product
of the symmetries about these leaves is loxodromic as long as
$\theta\in[-\pi/2,\pi/2]$ (see Lemma \ref{FalZoc}).  More precisely,
we will show that for these values of $\theta$, the isometry
associated to the matrix $M^\theta_{2,r_3}\overline{M^\theta_{3,r_3}}$
is loxodromic for any pair $(r_2,r_3)\in \R_{>0}^2$. Using the above
matrix form, it is seen that the trace of this matrix is

\begin{equation}
  \tr M^\theta_{2,r_2}\overline{M^\theta_{3,r_3}} = 
  2r_3^2e^{i\theta}+\frac{2}{r_2^2}e^{-i\theta}+1+r_2^2r_3^2+\frac{1}{r_2^2r_3^2}+\frac{r_3^2}{r_2^2}.
\end{equation}

\n This yields

\begin{eqnarray}
  \Re\left(\tr M^\theta_{2,r_2}\overline{M^\theta_{3,r_3}}\right) & = &
  2r_3^2\cos{\theta}+\frac{2}{r_2^2}\cos{\theta}+1+r_2^2r_3^2+\frac{1}{r_2^2r_3^2}+\frac{r_3^2}{r_2^2}\nonumber\\
  & \geq & 1+r_2^2r_3^2+\frac{1}{r_2^2r_3^2}\mbox{ while }\cos{\theta}\geq 0\nonumber\\
  & \geq & 3
\end{eqnarray}

\n This implies that the isometry associated to
$M^\theta_{2,r_3}\overline{M^\theta_{3,r_3}}$ is loxodromic as long as
$\theta\in[-\pi/2,\pi/2]$ and for any pair $(r_2,r_3)\in\R_{>0}^2$, as
shown by Remark \ref{trace3}. As a consequence of Lemma \ref{FalZoc},
the corresponding leaves of $S_2$ and $S_3$ are disjoint.
\end{proof}

\begin{rem}
  In the case of \PSL, or more generally in the case of real split Lie
  groups, it is possible to prove the discreteness of the image of
  $\rho$ by studying the coordinate changes induced by the flip moves:
  these moves preserve the positivity of cross-ratios, and this leads
  to the discreteness of $\rho$ (see for instance \cite{FockGon3}
  pages 87 to 89). Such an approach is not possible here. Notice for
  instance that if $(a,b,c)$ and $(c,d,a)$ are real ideal triangles
  sharing an edge, then the two ideal triangles obtained after a flip
  move, namely $(a,b,d)$ and $(b,c,d)$ are not real in general.
\end{rem}
\section{Remarks and comments\label{rem}}
\subsection{The case of real positive decorations : $\R$-Fuchsian representations.}
Let us focus for a moment on the special case where the bending
decoration is positive: for all edge $e$ of $T$, $\tD(e)\in\R>0$. It
this case, all the triangles constructed from $\tD$ are contained in
the standard real plane $\HdR$. As mentionned in Remark
\ref{linkclassical}, the $\tZ$-invariant is in this case the usual
cross-ratio in the upper half-plane. We recover this way the classical
\textit{shear coordinates}, and the action of the
$\rho(\pi_1(\Sigma))$ on the upper half-plane, when $\rho$ is a
discrete and faithful representation in $\PSL$. This corresponds to
the embedding $\PSL\sim\mbox{PO(2,1)}$ as the stabilizer of
$\HdR$. Moreover, when $z\in\R_{>0}$, the restriction of the real
symmetry $\sigma_z$ to $\HdR$ is a half-turn. We recover thus also the
explicit combinatorial description of classes of discrete and faithful
representations in \PSL\, given for instance by Fock and Goncharov in
\cite{FockGonA} by means of elementaries isometries (see section
\ref{embed}).  Note that the parabolicity criterion for peripheral
homotopy classes in \cite{PenR} or \cite{FockGonA} is the same as here
(it is expressed in a additive way in \cite{PenR}, where the situation
is slightly different, and the coordinates are expressed using
logarithms of cross-ratios).  In this particular case, if $\tau$ is
one of the real ideal triangles constructed from $\tD$, the prism
$\mathfrak{p}_\tau$ is the inverse image of $\tau$ by the orthogonal
projection onto $\HdR$.

\subsection{Embeddings of the Teichm\"uller in the PU(2,1)-representation variety \label{embed}}
Let us go back for a moment to the case of representations in \PSL,
the group of holomorphic isometries of the complex hyperbolic line
$\HuC$. In this frame, we can define a $\HuC$-realisation of the Farey
set of a cusped surface $\Sigma$ as a pair $(\phi,\rho)$, where
$\rho:\pi_1(\Sigma)\longrightarrow\PSL$ is a discrete and faithful
representation and $\phi$ is a $\rho$-equivariant mapping from the
Farey set to the boundary of the Poincar\'e disc. Denote by
$\mathcal{DF}$ the set of \PSL-classes of discrete and faithful
representations of $\pi_1(\Sigma)$ in \PSL, and by $\mathcal{DF}^+$
the set of \PSL-classes of $\HuC$-realizations of $\Fa(\Sigma)$.

Let $m$ be a point of $\Fa(\Sigma)$, corresponding to a fixed point of
a parabolic $c$, and let $(\phi,\rho)$ be a $\HuC$-realisation. The
$\rho$-equivariance of $\phi$ implies that $\phi(m)$ is fixed by
$\rho(c)$. Now, $\rho$ being discrete and faithful, $\rho(c)$ is
either parabolic or loxodromic. When $\rho(c)$ is hyperbolic,
$\phi(m)$ may be any of the two fixed points of $\rho(c)$. Consider
the projection

\begin{eqnarray}\bp & :
  &\mathcal{DF}^+\longrightarrow\mathcal{DF}\nonumber\\&
  &[(\phi,\rho)]\longmapsto[\rho]\label{projection}
\end{eqnarray}

\n Let $[[1,n]]$ be the set of integers between $1$ and
$n$. For any subset $I=\{i_1,\cdots,i_k\}$ of $[[1,n]]$, define

$$\mathcal{P}_I=\left\{[\phi,\rho]\in\mathcal{DF}^+|\rho(c_i) \mbox{ is parabolic
  }\Leftrightarrow i\in I\right\}.$$

\n Then $\mathcal{DF}^+$ decomposes as the disjoint union

\begin{equation}\label{stratif}
\mathcal{DF}^+=\coprod_{I\subset [[1,n]]}\mathcal{P}_I,
\end{equation}
and the restriction to $\mathcal{P}_I$ of the projection
(\ref{projection}) is $2^{n-|I|}$ to 1. In particular, it is $2^n$ to
1 when restricted to $\mathcal{P}_\emptyset$, which is the set of
realisations associated to totally hyperbolic representations, and it
a bijection when restricted to $\mathcal{P}_{[[1,n]]}$, which
corresponds to the Teichm\"uller space.

Once an ideal triangulation $T$ of $\Sigma$ is fixed, shear
coordinates provide a bijection between the set of positive decoration
of $T$ (that is, mappings $\td:e(T)\longrightarrow \R_{>0}$), and the
set of $\mathcal{DF}^+$. The main tool is the classical cross-ratio,
used as a gluing invariant of two ideal triangles in $\HuC$. It is
also possible to give an explicit representative for a representation
associated to a given decoration by use of elementary isometries. This
time, the elementary isometries are

$$I_x=\begin{bmatrix} 0 & \sqrt{x}\\ -1/\sqrt{x} & 0\end{bmatrix},$$
for an edge of type 1 intersecting an edge of $T$ decorated by the
positive number $x$, and
$$E=\begin{bmatrix} 1 & 1 \\ -1 & 0\end{bmatrix},$$
for a positively oriented edge of type 2.  The mechanic of the
construction is the same as what we did in section
\ref{sectionbending}, only simplified by the fact that both types of
elementary isometries are holomorphic, thus there is no need of
colouring faces of $\hat T$ in the classical case. This material is
classical and exposed for instance in \cite{FockGonA}. Notice that if
$c_j$ is a peripheral homotopy class around the deleted point $x_j$,
the parabolicity of $\rho(c_j)$ is equivalent to the condition that
the associated positive decoration is \textit{balanced at $x$} (that
is, the product of all positive numbers on edges adjacent to $x$
equals 1). Type-preserving representations, and therefore the
Teichm\"uller space of $\Sigma$ correspond to positive decorations
which are balanced at every deleted point of $\Sigma$. We call such
decorations simply \textit{balanced}.

Fix a bipartite ideal triangulation $T$. The set of positive decorations of $T$
is $\R_{>0}^{\sharp e(T)}$. To any real number
$\theta$ is associated a mapping 

\begin{eqnarray}\psi_\theta & : & \R_{>0}^{\sharp e(T)}  \longrightarrow \mathcal{BD}_T\nonumber\\
& &\td \longmapsto \tD=\td e^{i\theta}.\label{psitheta}
\end{eqnarray}
 
This mapping induces a mapping from $\mathcal{DF}^+$ to
$\mathcal{BR}_T$, which maps the realisation associated to $\td$ to
the $T$-bent realisation associated to the regular bending decoration
$\td e^{i\theta}$. Restricting this induced mapping to those
$\HuC$-realisation corresponding to balanced positive decorations, we
can rephrase Theorem \ref{theodiscret} as follows.

\begin{theo}\label{theoembed}
  Let $\theta\in[-\pi/2,\pi/2]$ be a real number and $T$ be a
  bipartite ideal triangulation of $\Sigma$.  The mapping
  $\psi_\theta$ defined in (\ref{psitheta}) induces a pair of
  embeddings of $\mathcal{DF}(\Sigma)$ of $\Sigma$ in Hom($\pi_1$,
  PU(2,1))/PU(2,1) of which images contain only classes of discrete
  and faithful representations. 
\end{theo}
\begin{proof}
  Restricting the mapping $\td\longmapsto \td e^{i\theta}$ to balanced
  decorations of $T$ produces discrete, faithful and type-preserving
  representations of $\pi_1(\Sigma)$ with images contained in PU(2,1)
  since $T$ is bipartite.  Once a coloring of the faces of $\hat T$
  is fixed, we obtain two injective applications by mapping the point
  in $\mathcal{T}(\Sigma)$ associated to $\td$ to the class of
  representations associated to $\td e^{i\theta}$ corresponding either
  to white triangles or to black triangles. These two embeddings are
  identified by the complex conjugation in $\HdC$, and correspond in
  fact to a single embedding in Hom($\pi_1$,PU(2,1))/Isom($\HdC$).
\end{proof}

Note that Theorem \ref{theodiscret} states as well that the
parabolicity of images of peripheral loops is preserved by
$\psi_\theta$, and thus the image of $\psi_\theta$ admits a similar
decomposition as (\ref{stratif}).

\subsection{Link with previously known families of examples.}
In this section, we draw the connection between $T$-bent realizations
and families of examples described in the previous works
\cite{FK1,GuP2,Wi2}.

\paragraph{The 1-punctured torus.} 

\begin{figure}  
\begin{center}
\begin{pspicture}(-2,-2)(4,4)
\includegraphics[width=5cm,height=5cm]{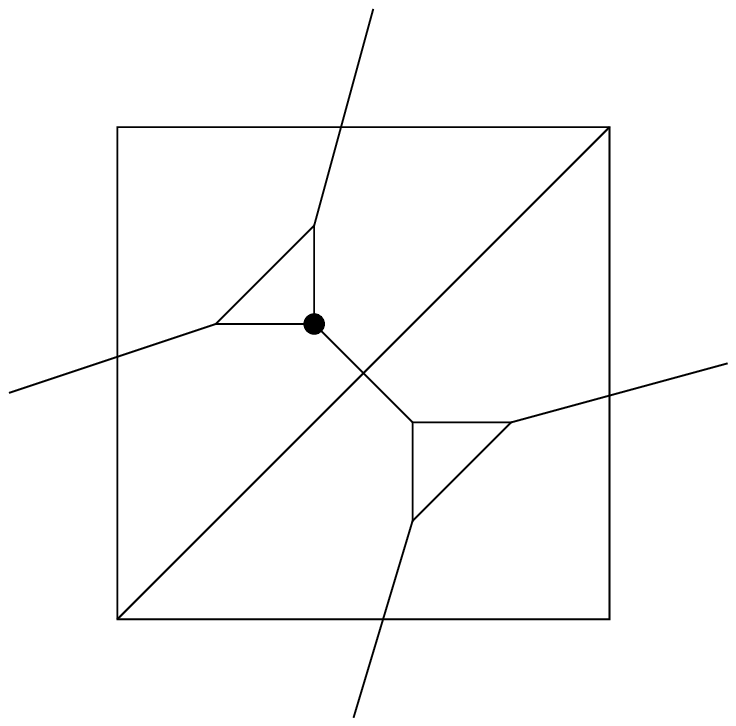}
\rput(-1,0.9){{\tt b}}
\rput(-4.1,4){{\tt w}}
\rput(-3,2.6){v}
\rput(-4,1.9){$e_1$}
\rput(-3.4,0.9){$e_3$}
\rput(-2.35,3){$e_2$}
\psarc{->}(-1.2,3.5){1}{0}{60}
\psline{->}(-3.7,2.3)(-1,2.3)
\psline{->}(-1.5,1.2)(-1.5,4)
\end{pspicture}
\vspace{-2cm}
\caption{The 1-punctured torus\label{pictor}}
\end{center}
\end{figure}

In this case $T$ consists of two triangles, as indicated on figure
\ref{pictor}. We will use the vertex $v$ marked on the figure as
basepoint. There are two faces, of which colour is indicated by ${\tt
  w}$ and ${\tt b}$ on figure \ref{pictor}, and three edges, labelled
by $e_1$, $e_2$ and $e_3$ on figure \ref{pictor}. In the case of a
regular bending decorations,the decoration is given three positive
real numbers $x_1$, $x_2$ and $x_3$ and $\theta\in[0,2\pi[$ such that
the edge $e_i$ is decorated by $x_i,\theta$. Following the results of
section \ref{explicit}, we see that the identifications between
opposite faces of the square correspond to the following holomorphic
isometries of $\HdC$. Call $A$ and $B$ the isometries associated
respectively to the horizontal and vertical identifications of the
opposite sides of the square. Following section \ref{explicit}, these
isometries are given by

\begin{equation}
  \left\{\begin{matrix}\label{AB}
      A & = &\mathcal{E}\circ\sigma_{x_2,\theta}\circ \mathcal{E}^{-1}\circ\sigma_{x_2,\theta}\\
      \\
      B & = &\sigma_{x_2,\theta}\circ \mathcal{E}^{-1}\circ\sigma_{x_3,\theta}\circ \mathcal{E}.
\end{matrix}\right.
\end{equation}

\n As a consequence, we see that the group $\la A,B\ra$ has index two
in the group generated by the three real symmetries
$I_1=\mathcal{E}\circ\sigma_{x_1,\theta}\circ \mathcal{E}^{-1}$,
$I_2=\sigma_{x_2,\theta}$ and $I_3=
\mathcal{E}^{-1}\circ\sigma_{x_3,\theta}\circ \mathcal{E}$. The group
$\la I_1,I_2,I_3\ra$ is an example of a so-called Lagrangian triangle
group. This example of bending has been exposed with a different point
of view in \cite{Wi2} (see also\cite{Wi6}).
 
In \cite{Wi2}, the discreteness result is stated with an angle
$\alpha\in[-\pi/4,\pi/4]$.  This angle $\alpha$ is actually half the
bending parameter $\theta$ we use here. It may be interpreted as an
angle between a real ideal triangle $\Delta$ and the splitting surface
Spl($\Delta,\Delta'$), where $\Delta'$ is adjacent to $\Delta$. From
this point of view, Spl($\Delta$,$\Delta'$) is bisecting the pair
$(\Delta,\Delta')$.
 
\paragraph{The Toledo invariant and the examples of Gusevskii and Parker}
The Toledo invariant is a conjugacy invariant defined for
representations of fundamental groups of closed surfaces, and for
type-preserving representations of cusped surfaces. We refer the
reader to \cite{Tol} and \cite{KoMa,GuP2} for its definition and main
properties. Let us just recall that if $\rho$ is such a
representation, then

\begin{itemize}
\item if $\Sigma$ has punctures, then ${\bf tol}(\rho)$ is a real number in
  the interval $[-\chi,\chi]$, where $\chi$ is the Euler
  characteristic of $\Sigma$,
\item if not, then ${{\bf tol}}(\rho)$ belongs to $2/3\mathbb{Z}\cap
  [\chi,-\chi]$.
\end{itemize}

Let $(\phi,\rho)$ be a $T$-bent realization of $\mathcal{F}_\infty$,
where $T$ is a bipartite triangulation, and $\Omega$ be a fundamental
domain for the action of $\pi_1(\Sigma)$ on $\tilde\Sigma$. We might
see $\Omega$ as a family of triangles $(\Delta_1,\cdots,\Delta_m)$.
Then it follows from \cite{GuP2,Tol} that the Toledo invariant 
${\bf tol}(\rho)$ equals twice the sum of the Cartan invariants of the ideal
triangles $\phi(\Delta_i)$. In our particular case, all the triangles
are real. We obtain therefore directly the

\begin{prop}\label{Toledo}
  Let $(\phi,\rho)$ be a $T$-bent realization of $\mathcal{F}_\infty$,
  with $\rho$ type-preserving.  The Toledo invariant of $\rho$ is
  equal to zero.
\end{prop}

In \cite{GuP2}, Gusevskii and Parker have described for each genus $g$
and number of punctures $n$ a 1-parameter family $(\rho_t)_{t\in
  [-\chi,\chi]}$ of non PU(2,1)-equivalent discrete, faithful and
type-preserving representations of a Riemann surface of genus $g$ with
$n$ punctures having the property that the Toledo invariant of
$\rho_t$ equals $t$.  This shows that all the possible values of the
Toledo invariant for non-compact surfaces are realised by discrete
and faithful representation. To prove this result, Gusevskii and
Parker start from discrete and faithful representations of the modular
group in PU(2,1) and pass to a finite index subgroup using
Millington's theorem (see \cite{GuP2}). In their construction, they
show that $\rho_0$ preserves a real plane (this is a so-called
$\R$-Fuchsian representation).  Therefore $\rho_0$ is the unique
intersection between Gusevskii and Parker's family of representations
and our one.

\paragraph{The 3-punctured sphere and the examples of Falbel and Koseleff.}

\begin{figure}  
\begin{center}
\begin{pspicture}(-2,-2)(4,4)
\includegraphics[width=5cm,height=5cm]{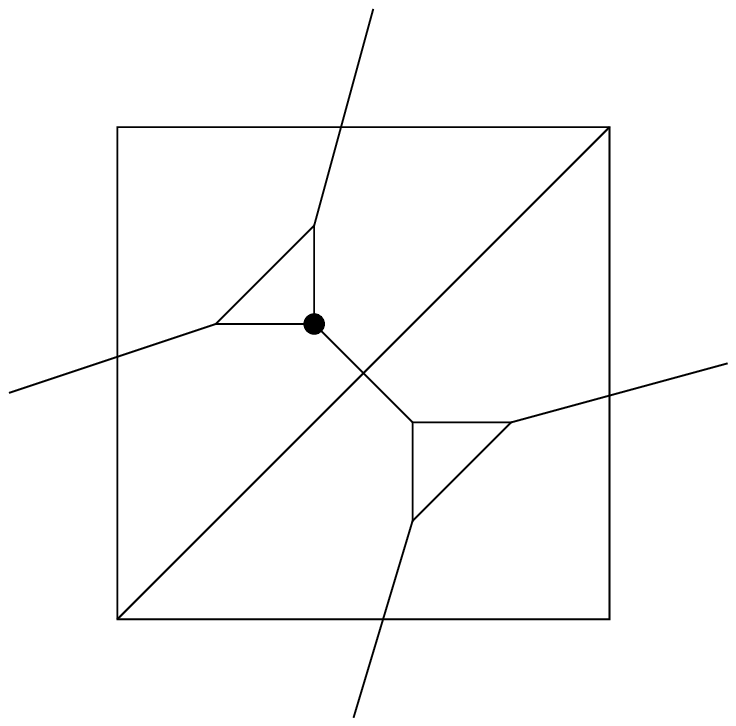}
\rput(-1.2,3.5){{\tt b}}
\rput(-3.5,1.7){{\tt w}}
\rput(-3,2.5){v}
\rput(-4,1.9){$e_3$}
\rput(-3.4,0.9){$e_3$}
\rput(-2.35,3){$e_2$}
\rput(-0.5,3){$e_1$}
\rput(-1.4,4.3){$e_1$}
\psarc{->}(-1.2,3.5){1}{0}{60}
\end{pspicture}
\vspace{-2cm}
\caption{The 3-punctured sphere\label{picsphere}}
\end{center}
\end{figure}

This time we are using the bipartite triangulation of the 3-punctured
sphere showed on figure \ref{picsphere}. The representation of the
fundamental group associated to the decoration given by
$\delta(e_i)=x_i$ and $\alpha(e_i)=\theta_i$ is given by

\begin{equation}
  \left\{\begin{matrix}\label{ABC}
      A & = &\mathcal{E}^{-1}\circ\sigma_{x_1,\theta_1}\circ \mathcal{E}^{-1}\circ\sigma_{x_2,\theta_2}\\
      \\
      B & = &\sigma_{x_2,\theta_2}\circ \mathcal{E}^{-1}\circ\sigma_{x_3,\theta_3}\circ \mathcal{E}^{-1}.\\
      \\
      C& = & \mathcal{E}\circ\sigma_{x_3,\theta_3}\circ \mathcal{E}^{-1}\circ \sigma_{x_1,\theta_1}\circ \mathcal{E},
\end{matrix}\right.
\end{equation}

\n It is easily checked that $ABC=1$. Using the matrices given in
section \ref{explicit}, we see that the representation is type
preserving if and only if $x_1=x_2=x_3=1$ and none of the $\theta_i$'s
is equal to $\pi$. When $\theta_1=\theta_2=\theta_3\in[-\pi/2,\pi/2]$,
this provides through theorem \ref{theodiscret} a 1-parameter family
of discrete, faithful and type-preserving representations of the
fundamental group of the 3-punctured sphere.

\n Moreover, it is possible to prove that in the case where
$\delta(e_i)=1$ and $\alpha(e_i)=\theta$ for all $i$, then there
exists three real symmetries $s_1$, $s_2$ and $s_3$ such that
$A=s_1s_2$ and $B=s_2s_3$. Call $Q_i$ the mirror of $s_i$. Since $A$
and $B$ are parabolic, the mirrors of the $s_i$'s are mutually
asymptotic, that is $Q_i\cap Q_{i+1}$ consists of exactly one point in
$\partial\HdC$. Therefore these groups belong to the family of groups
studied by Falbel and Koseleff in \cite{FK1}. Note moreover that the
discreteness of these groups was not proved in \cite{FK1}, where the
focus is on deformations of groups preserving a complex line.

\end{document}